\documentclass[12pt]{article}
\usepackage{amsmath}
\usepackage{amssymb}
\usepackage{amsthm}
\usepackage{latexsym}
\usepackage{fancyhdr,lastpage,graphicx}
\theoremstyle{plain}
\newtheorem{thm}{Theorem}[section]
\newtheorem{lemma}{Lemma}[section]

\newtheorem{prop}[lemma]{Proposition}

\theoremstyle{remark}
\newtheorem{rem}[lemma]{Remark}

\newtheorem{example}[lemma]{Example}
\def\PP{\mathbb{P}}
\def\ve{{\varepsilon}}
\begin{document}
\thispagestyle{plain}
\begin{center}

\vspace{0.2in}

\LARGE \textbf{Structure of Large Random Hypergraphs}

\vspace{0.2in}

\large {\bfseries R.W.R. Darling
\footnote{National Security Agency, P.O. Box 535, Annapolis Junction, Maryland, 20701-0535, USA. Email: rwrd@afterlife.ncsc.mil} 
and J.R. Norris
\footnote{Statistical Laboratory, Centre for Mathematical Sciences, Wilberforce Road,Cambridge, CB3 0WB, UK. Email: j.r.norris@statslab.cam.ac.uk}}

\vspace{0.2in}
\small
\today

\end{center}

\vspace{0.2in}

\begin{abstract}
The theme of this paper is the derivation of analytic formulae
for certain large combinatorial structures.
The formulae are obtained via fluid limits
of pure jump type Markov processes, established under simple conditions
on the Laplace transforms of their L\'evy kernels. 
Furthermore, a related Gaussian approximation allows us to describe 
the randomness which may persist in the limit when certain parameters 
take critical values. 
Our method is quite general,
but is applied here to vertex identifiability in random hypergraphs.
A vertex $v$ is identifiable in $n$ steps if there is a hyperedge
containing $v$ all of whose other vertices are identifiable in fewer
steps. We say that a hyperedge is identifiable if every one of its vertices is
identifiable.
Our analytic formulae describe the asymptotics of the number of identifiable vertices 
and the number of identifiable hyperedges for a Poisson($\beta$) random hypergraph
$\Lambda$ on a set $V$ of $N$ vertices, in the limit as $N\rightarrow\infty$.
Here  $\beta$ is a formal power series with non-negative coefficients 
$\beta_{0}, \beta_{1}, \ldots$, and $(\Lambda (A))_{A \subseteq V}$ 
are independent Poisson random variables such that 
$\Lambda (A)$, the number of hyperedges on $A$, 
has mean $N\beta_j /\binom{N}{j}$ whenever 
$\left| A \right| = j$.

\emph{Keywords} hypergraph, component, cluster, Markov process, random graph

\emph{AMS (2000) Mathematics Subject Classification}: Primary 05C65; Scondary  60J75, 05C80

\end{abstract}
\section{Introduction}
\subsection{Motivation}
We are interested in the evolution of certain statistically symmetric random
structures, extended over a large finite set of points, when points are
progressively removed in a way which depends on the structure. The
initial condition of the structure may allow few possibilities for the
removal of points, indeed it may be that, once a small proportion of
points are removed, the process terminates. On the other hand, the
removal of points may cause the structure to ripen, eventually yielding
a large proportion of the initial points. Our analysis will enable us to
demonstrate a sharp transition between these two sorts of behaviour as
certain parameters pass through critical values.

Let us illustrate this phenomenon by a simple special case. Consider
the complete graph on $N$ vertices and declare each vertex to be open
with probability $p$, each edge to be open with probability $\alpha/N$.
Suppose that we are allowed to select an open vertex, remove it, and 
declare open any other vertices sharing an open edge with the selected vertex.
If we continue in this way until no open vertices remain, we eventually
remove every vertex connected to an open
vertex by open edges. We shall see that the proportion of vertices thus
removed converges in probability as $N\to\infty$ and that
the limit $z^*(p,\alpha)$ is the unique root in $[0,1)$ of the equation
\begin{equation*}
\alpha  z+\log(1-z)=\log(1-p).
\end{equation*}
Thus, for small values of $p$, there is a dramatic change in behaviour
as $\alpha$ passes through 1. As $p\downarrow0$, for $\alpha \le 1$,
\begin{equation*}
z^*(p,\alpha)/p\to1/(1-\alpha)
\end{equation*}
but for $\alpha>1$
\begin{equation*}
z^*(0+,\alpha)>0.
\end{equation*}
Of course this is a reflection of well known connectivity properties of
random graphs, discovered by Erd{\H{o}}s and R{\'e}nyi \cite{ER}, and
discussed, for example, in \cite{B84}.

The class of models considered in this paper is a natural
generalization of some classical models of random graphs and hypergraphs,
which may be further motivated as follows.
Phase transitions in combinatorial problems constitute an area of active research among computer scientists. Many ``hard'' combinatorial problems can be cast as \emph{satisfiability} problems, which seek to assign a truth value to each of a set of Boolean variables, such that a collection of logical conjunctions are simultaneously satisfied. 
Phase transitions for random satisfiability (``random k-SAT'') problems have 
been studied by researchers at Microsoft \cite{A99}, \cite{AKKK99}, \cite{W02} 
and IBM \cite{CGHS03}, but difficult questions remain unanswered. 
The random hypergraph model herein may be viewed as a simplification of the random satisfiability model: a vertex corresponds to a Boolean variable, and a hyperedge to 
the set of variables appearing in a specific logical conjunction, 
neglecting the truth or falsehood assigned to those variables. Under this 
simplification, definitive critical parameters are obtained which 
shed light on the random satisfiability model, and whose derivation may serve 
as a template for analysis of mixed
satisfiability problems.  
\subsection{Hypergraphs}
Let $V$ be a finite set of $N$ vertices. By a {\em hypergraph} on $V$ we
mean any map
\begin{equation*}
\Lambda:\mathcal{P}(V)\to\mathbb{Z}^+.
\end{equation*}
Here $\mathbb{Z}^+$ denotes the set of non-negative integers.
The reader may consult \cite{ENC} for an overview of the theory of hypergraphs:
however the direction pursued here is largely independent of previous work.
We emphasise that, in distinction to much of the combinatorial literature
on hypergraphs, we allow the possibility that more than one edge is assigned
to a given subset, thus we are considering multi-hypergraphs. Moreover we
do not insist that all hyperedges have the same number number of vertices.
Much of the literature is restricted to this {\em uniform} case. Our
methods allow a significant broadening of the class of models for which
asymptotic computations are feasible.
Hyperedges over vertices are called {\em patches\/}
(loops in \cite{ENC}) and hyperedges over $\emptyset$ are called {\em debris\/}.
The total number of hyperedges is
\begin{equation*}
|\Lambda|=\sum_A\Lambda(A).
\end{equation*}

\subsection{Accessibility and Identifiability}
Interest in large random graphs has often focused on the sizes of their
connected components. If there is given also, as in the example above,
a set of distinguished vertices $V_0$, then it is natural to seek to determine
the proportion of all vertices connected to $V_0$. 

In the more general context of hypergraphs there is more than one interesting
counterpart of connectivity. Given a hypergraph $\Lambda$ on a set $V$, 
we say that a vertex $v$ is {\em accessible in $1$ step} or, equivalently,
{\em identifiable in $1$ step} if $\Lambda(\{v\})\ge1$.
We say, for $n=2,3,\dots$, that a vertex is {\em accessible in $n$ steps} if 
it belongs to some subset $A$ with $\Lambda(A)\ge1$, 
{\em some} other element of which is accessible in less than $n$ steps.
A vertex is {\em accessible} if it is accessible in $n$ steps for some $n\ge1$.

On the other hand, we say that a vertex is {\em identifiable in $n$ steps} if 
it belongs to some subset $A$ with $\Lambda(A)\ge1$, 
{\em all} of whose other elements are identifiable in less than $n$ steps.
A vertex is {\em identifiable} if it is identifiable in $n$ steps for some $n\ge1$.

The notion of accessibility may be appropriate to some physical models similar to
percolation, whereas identifiability is more relevant to knowledge-based structures.
We shall examine only the notion of identifiability.

Given a hypergraph $\Lambda$ without patches and a distinguished vertex $v_0$,
we say that a vertex $v$ is {\em accessible from $v_0$} if it is accessible in the
hypergraph $\Lambda+1_{\{\{v_0\}\}}$, that is, in the hypergraph obtained
from $\Lambda$ by adding a single patch at $v_0$. 
Identifiability from $v_0$ is defined
similarly. The set 
of vertices accessible from $v_0$ is the {\em component} of $v_0$, as studied in
\cite{ENC,KL93,K85,SPS85}. The set of vertices identifiable from $v_0$ 
is the {\em domain} of $v_0$, as studied by Levin and the current authors 
\cite{DL01}. 
We shall not consider further in this paper these vertex-based notions.

The process of identification is dual to the process leading to the 2-core
of a graph or hypergraph, that is to say, the maximal subgraph in which every
non-isolated vertex has degree at least 2. 
In the former process one removes vertices having a
1-hyperedge, in the latter one removes edges containing
a vertex of degree one. In this duality, non-identifiable vertices
correspond to the 2-core. Thus our results may be interpreted as giving
the asymptotic size of the 2-core for a certain class of random hypergraphs.

\subsection{Hypergraph Collapse}
\begin{figure}
\begin{center}
\scalebox{0.6}{\includegraphics{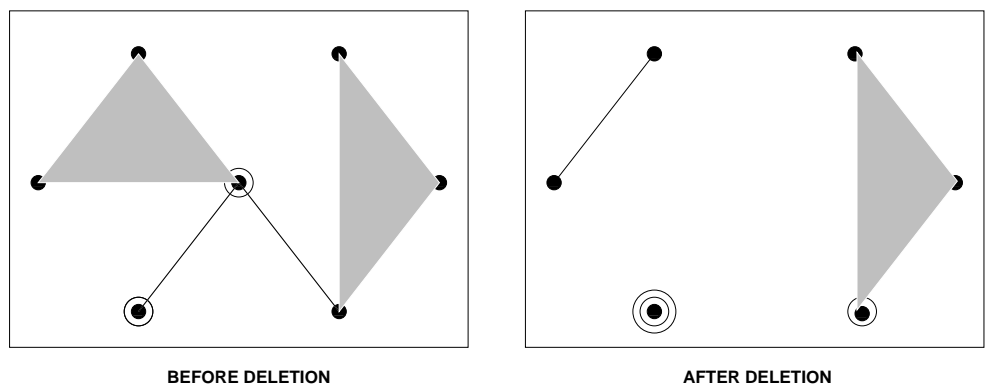}}
\caption{EXAMPLE OF A PERMITTED COLLAPSE - DELETION OF ONE VERTEX}\label{f1}
\end{center}
\end{figure}

It will be helpful to think of the identification  of vertices as a
progressive activity. Once a vertex is identified, it is
\emph{removed} or \emph{deleted} from the vertex set, in a manner which is explained
below. Thus, we shall consider an evolution of hypergraphs by the removal of
vertices over which there is a patch. A hypergraph with no patches
will therefore be {\em stable\/}. Given a hypergraph $\Lambda$ and a vertex $v$,
we can arrive at a new hypergraph $\Lambda'$ by removing $v$ from each
of the hyperedges of $\Lambda$. Thus
\begin{equation*}
\Lambda'(A)=
\begin{cases}
\Lambda(A)+\Lambda(A\cup\{v\}) & \text{if $v\notin A,$}\\
0 & \text{if $v\in A$}.
\end{cases}
\end{equation*}
For example, in Figure \ref{f1}, the patch on the central vertex is selected,
and that vertex is removed; this causes a triangular face
to collapse to an edge, and two edges incident to the vertex to
collapse to patches on the vertices at the other ends. Note that this leaves
 two patches on the lower left vertex.

%
If $\Lambda(\{v\})\ge 1$ then we say that $\Lambda'$ is obtained
from $\Lambda$ by a ({\em permitted\/}) {\em collapse\/}. Starting from
$\Lambda$, we can obtain, by a finite sequence of collapses, a stable
hypergraph $\Lambda_\infty$. Denote by $V^*$ the set of vertices removed
in passing from $\Lambda$ to $\Lambda_\infty$. The elements of $V^*$ 
are the identifiable vertices. We write $\Lambda^*$ for the 
{\em identifiable hypergraph}, given by
\begin{equation*}
\Lambda^*(A)=\Lambda(A)1_{A\subseteq V^{*}}.
\end{equation*}  
We note that $V^*$, and hence $\Lambda^*$ and $\Lambda_\infty$, do not
depend on the particular sequence of collapses chosen. For, if $v_1,
v_2,\dots$ and $v_1', v_2',\dots$ are two such sequences, and if
$v_n\neq v_k'$ for all $k$, then we can take $n$ minimal and find $k$ such
that $\{v_1,\dots, v_{n-1}\}\subseteq\{v_1', \dots,v_k'\}$; then, with
an obvious notation, $\Lambda_k'(\{v_n\})\ge
\Lambda_{n-1}(\{v_n\})\ge 1$, so $v_n$ must, after all, appear in
the terminating sequence $v_1', v_2',\dots$, a contradiction. We note
also that $V^*$ increases with $\Lambda$.

\subsection{Purpose of This Paper}
The main question we shall address is to determine the asymptotic sizes
of $V^*$ and $\Lambda^*$ for certain generic random hypergraphs, as the
number of vertices becomes large.  We note that, since the number of
hyperedges is conserved in each collapse, all the identifiable hyperedges
eventually turn to debris:
\begin{equation*}
\Lambda_\infty (\emptyset) = |\Lambda^*|.
\end{equation*}

Note that $V^*$ depends only on $\min\{\Lambda, 1\}$.  In the case where
$\Lambda(A) =0$ for $|A| \ge 3$, the hypergraph min$\{\Lambda,
1\}$ may be considered as a graph on $V$ equipped with a number of
distinguished vertices.  Then $V^*$ is precisely the set of vertices
connected in the graph to one of these distinguished vertices.

\subsection{Poisson Random Hypergraphs}
Let $(\Omega, \mathcal{F}, \mathbb{P})$ be a probability space.  A {\em random
hypergraph\/} on $V$ is a measurable map
\begin{equation*}
\Lambda :\Omega \times \mathcal{P}(V)\to \mathbb{Z}^+.
\end{equation*}
An introduction to random hypergraphs may be found in \cite{KL93}, though we
shall pursue rather different questions here.
 We shall consider a class of random hypergraphs whose distribution is
determined by a sequence $\beta = (\beta_j : j \in \mathbb{Z}^+)$ of
non-negative parameters.  Say that a random hypergraph $\Lambda$ on $V$
is {\em Poisson\/}$(\beta)$ if
\begin{itemize}
\item The random variables $\Lambda(A), A \subseteq V$, are independent,
\item The distribution of $\Lambda(A)$ depends only on $|A|$,
\item $\sum_{|A|=j}  \Lambda(A) \sim$ Poisson$(N\beta_j)$,
$j=0,1,\dots, N$.
\end{itemize}
A consequence of these assumptions is that $\Lambda (A)$ has 
mean $N\beta_j /\binom{N}{j}$ whenever 
$\left| A \right| = j$.
Note that, when $N$ is large, for $j\ge2$, only a small fraction of the subsets
of size $j$ have any hyperedges, and those that do usually have just $1$. Also the
ratio of $j$-edges to vertices tends to $\beta_j$. Our assumption of Poisson
distributions is a convenient exact framework reflecting behaviour which holds
asymptotically as $N\rightarrow\infty$ under more generic conditions.

\subsection{Generating Function}
A key role is played by the power series
\begin{equation}\label{beta}
\beta(t) = \sum_{j \ge 0} \beta_j t^j
\end{equation}
and by the derived series
\begin{equation*}
\beta' (t)=\sum_{j \ge 1} j \beta_j t^{j-1},\quad \beta'' (t) =
\sum_{j\ge 2} j(j-1)\beta_j \ t^{j-2}.
\end{equation*}
Let $\beta$ have radius of convergence $R$.  The function $\beta'(t) +
\log (1-t)$ may have zeros in $[0,1)$ but these can accumulate only at
1.  Set
\begin{equation}\label{zstar}
z^* = \inf \{t \in[0,1): \beta'(t) + \log (1-t)<0\} \wedge 1
\end{equation}
and denote by $\zeta$ the set of zeros of $\beta'(t)+\log (1-t)$ in
$[0,z^*)$.  Note that if $\beta$ is a polynomial, or indeed if $R>1$,
then $z^*<1$.  Also, the generic and simplest case is where $\zeta$ is empty.
\section{Results}
We state our principal result first in the generic case.
\subsection{Hypergraph Collapse - Generic Case}
\begin{thm}
\label{t1}  
Assume that $z^* <1$ and $\zeta=\emptyset$.  For $N\in  \mathbb{N}$, let $V^N$
be a set of $N$ vertices and let $\Lambda^N$ be a Poisson$(\beta)$
hypergraph on $V^N$.  Then, as $N\to \infty$, the numbers of 
identifiable vertices and identifiable hyperedges satisfy the
following limits in probability :
\begin{equation*}
|V^{N*}|/N\to z^*,\quad |\Lambda^{N*}|/N\to \beta(z^*)-(1-z^*)\log (1-z^*).
\end{equation*}
\end{thm}
\begin{example}
The random graph with distinguished vertices, described in the
introduction corresponds to a Poisson$(\beta^N)$ hypergraph
$\Lambda^N$, where
\begin{equation*}
1-e^{-\beta^{N}_{1}}=p, \quad 1-e^{-2\beta^{N}_{2}/(N-1)}=\alpha/N
\end{equation*}
and $\beta^N_j=0$ for $j\ge 3$.  Note that  $\beta^N_1=\beta_1$
and $\beta^N_2\to \beta_2$ as $N\to\infty$, where $\beta_1 = -\log
(1-p)$ and $\beta_2 = \alpha/2$.  Theorem \ref{t1} extends easily to cases
where $\beta$ depends on $N$ in such a mild way:  one just has to check
that Lemma \ref{l71} remains valid and note that this is the only place that
$\beta$ enters the calculations.  We have $\beta(t) =-t\log
(1-p)+t^2\alpha/2$ so
\begin{equation*}
\beta'(t) + \log (1-t) = - \log (1-p) + t\alpha + \log (1-t).
\end{equation*}
Then $z^*$ is the unique $t\in [0,1)$ such that
\begin{equation*}
\alpha t + \log (1-t) = \log (1-p)
\end{equation*}
and $\zeta$ is empty, so $|V^{*N}|/N\to z^*$ in probability as $N\to
\infty$, as stated above.
\end{example}
\begin{example}
To illustrate critical phenomena, let $\beta(t) = \alpha(0.1 + 0.9t)^{7}$.
Let $x$, $y$, and $z$ refer to the re-scaled number of vertices eliminated, 
the number of patches, and the amount of debris, respectively;
 here ``re-scaled'' means after division by the number of vertices. 
Plots of $y$ and $z$ versus $x$  are shown in Figure \ref{f2}, 
for the choices $\alpha = 1185$ (solid) and $\alpha = 1200$ (dashed). 
In the case $\alpha = 1185$, $y$ hits zero when $x\approx 0.02$ , and so $z$
remains stuck at about 0.02. A very small increase in  $\alpha$, from $1185$
to $1200$, causes a dramatic change in the outcome:
after narrowly avoiding extinction (Figure \ref{f2}), the number of patches 
explodes (Figure \ref{f3}) as $x$ increases towards 1. 

Consider what the Figures tell us about the supercritical case 
$\alpha = 1200$: during the first 4\% of patch selections, 
there is rarely any other patch covering the same vertex as the one selected;
Figure \ref{f3} shows that, during the last 10\%
of patch selections, an average of $5792$ other patches
cover the same vertex as the one selected. [Read the labels 
on the $x$-axes carefully: Figure \ref{f2} is a close-up of the left-most
 4\% of the scale of Figure \ref{f3}.]
%
\begin{figure}[ht]
	\begin{center}
	\begin{tabular}{ccc}
	\scalebox{0.4}{\includegraphics{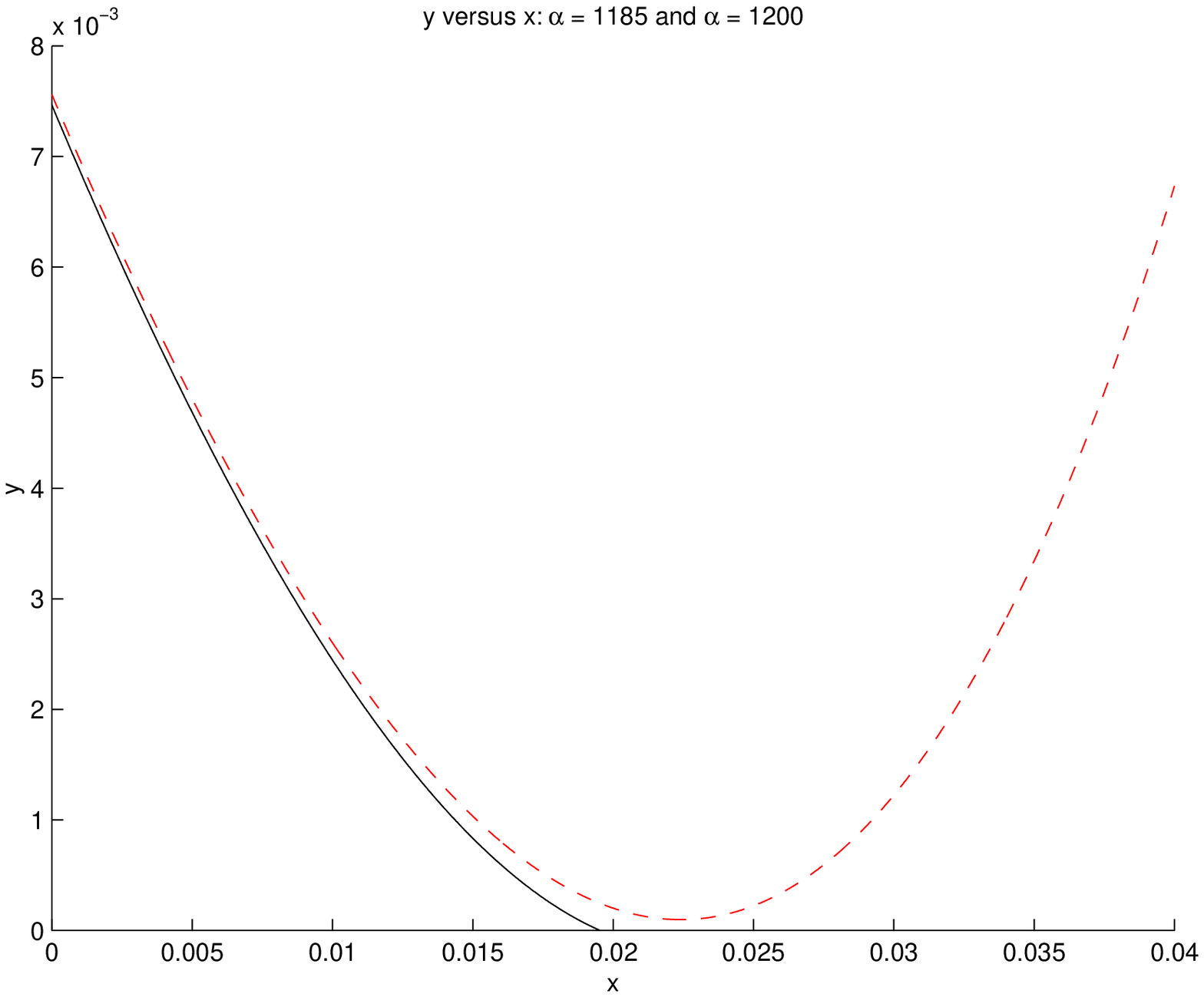}}
& \hspace{0.1in} &\scalebox{0.4}{\includegraphics{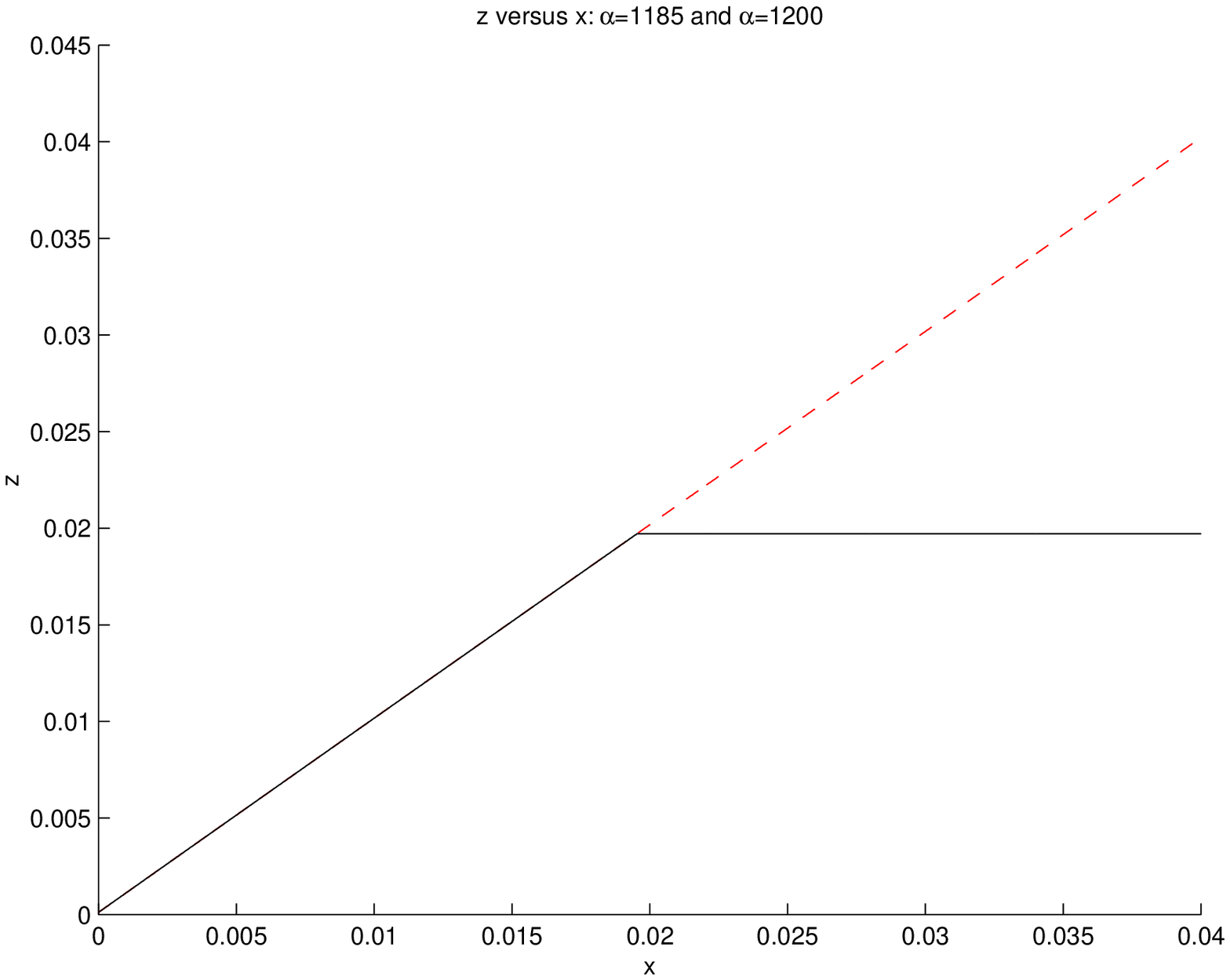}} \\
	\end{tabular}
	\end{center}
\caption{CHANGES IN BEHAVIOR NEAR A CRITICAL PARAMETER VALUE}\label{f2}
\end{figure}
\begin{figure}[ht]
	\begin{center}
	\begin{tabular}{ccc}
	\scalebox{0.4}{\includegraphics{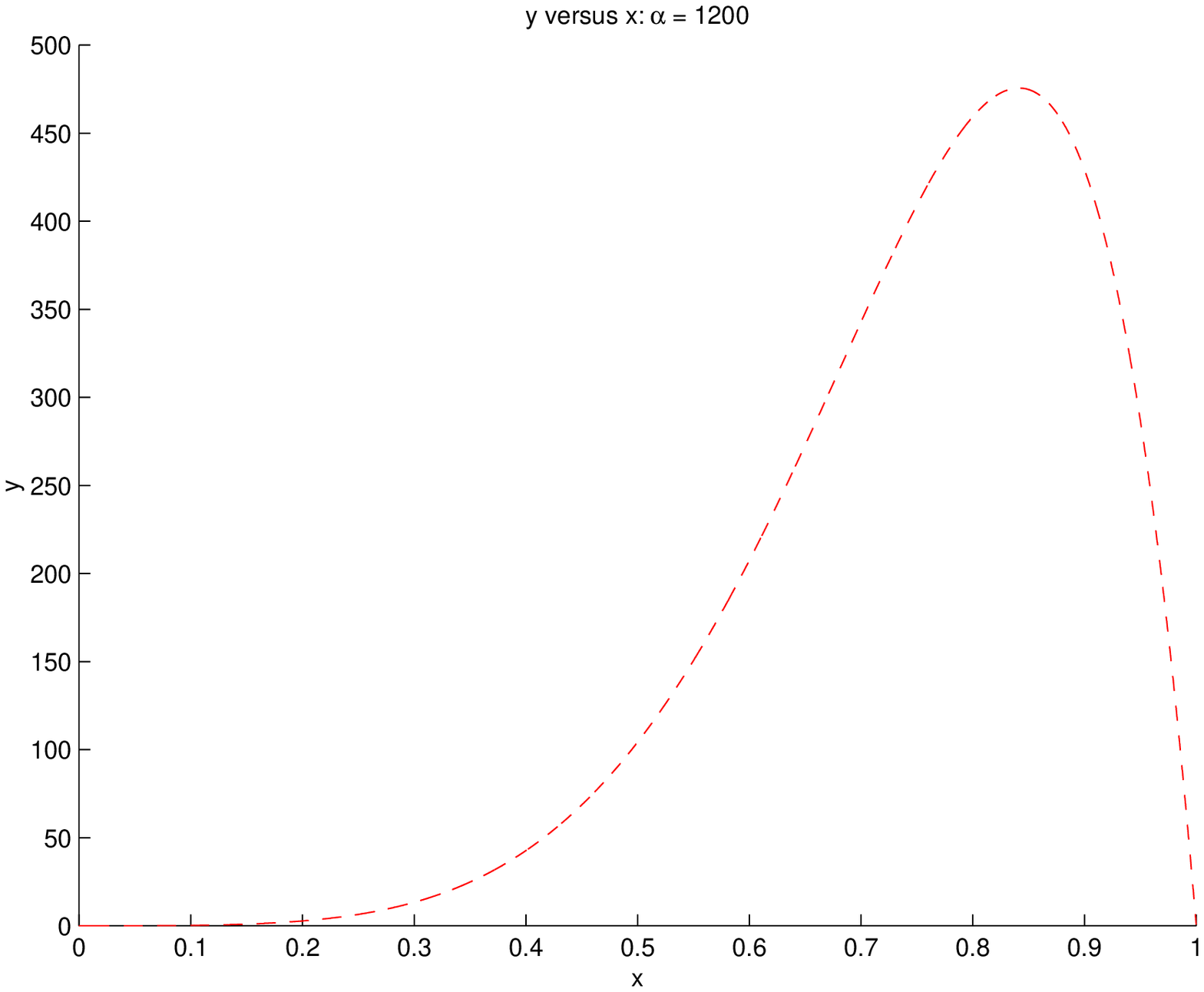}}
& \hspace{0.1in} &\scalebox{0.4}{\includegraphics{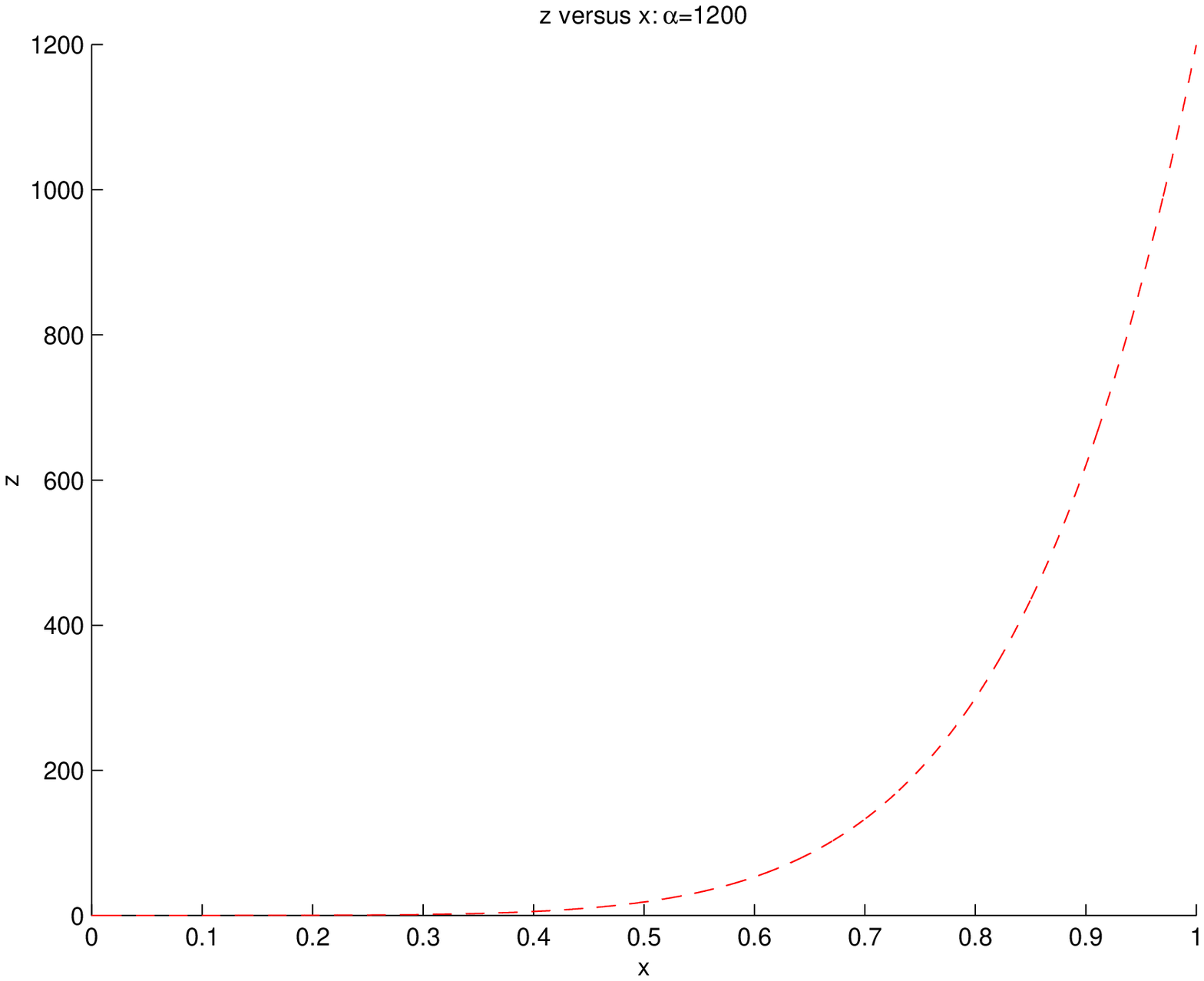}} \\
	\end{tabular}
	\end{center}
\caption{PATCHES AND DEBRIS IN THE SUPERCRITICAL REGIME}\label{f3}
\end{figure}
\end{example}
\subsection{Hypergraph Collapse - General Case}
In order to describe an extension of Theorem \ref{t1} to the case where
$\zeta$ is non-empty, we introduce the random variable
\begin{equation*}
Z=\min \{z\in \zeta:W\left(z/(1-z)\right)<0\} \wedge z^*
\end{equation*}
where $(W_t)_{t \ge 0}$ is a Brownian motion.

\begin{thm}\label{t2}
Assume that $R\notin \zeta$.  Then, for $V^{N*}$ and $\Lambda^{N*}$ as
in Theorem \ref{t1}, the following limits exist in distribution:
\begin{equation*}
|V^{N*}|/N\to Z,\quad |\Lambda ^{N*}|/N\to \beta(Z)-(1-Z)\log (1-Z).
\end{equation*}
\end{thm}
In the case where $\zeta$ has only a single point $\zeta_{0} < z^{*}$, then $Z$ is equal to 
 $\zeta_{0}$ with probability $\frac{1}{2}$ and equal to $z^{*}$  
with probability $\frac{1}{2}$.
We do not know what happens when $R\in \zeta$.  Proofs will be given in 
Section \ref{s7}.
\section{Randomized Collapse}\label{s4}
We introduce here a particular random rule for choosing the sequence of moves
by which a hypergraph is collapsed, which has the desirable feature that
certain key statistics of the evolving hypergraph behave as Markov chains.
It is by analysis of the asymptotics of these Markov chains as $N\to \infty$
that we are able to prove our main results.
\subsection{Induced Hypergraph}
Let $\Lambda$ be a Poisson$(\beta)$ hypergraph. 
For $S\subseteq V$ with $|S|=n$, let $\Lambda^S$ be the hypergraph
obtained from $\Lambda_0$ by removing all vertices in $S$.  Thus, for $A
\subseteq V\backslash S$ with $|A|=j$,
\begin{equation*}
\Lambda^S(A) = \sum_{B\supseteq A,\medspace B\backslash S=A} \Lambda (B) \sim
P\left(\lambda_j (N,n)\right)
\end{equation*}
where the Poisson parameter $\lambda_{j}(N,n)$ is computed as follows:
there are 
$\binom{n}{i}$ ways to choose $S\cap B$ such that $|B|=i+j$,
 and the Poisson parameter of $\Lambda_{0}(B)$ is
 $N \beta_{j+i}/ \binom{N}{i+j}$, so
\begin{equation*}
\lambda_{j}(N,n) = N\sum^n_{i=0} \beta_{j+i} \binom{n}{i}/\binom{N}{i+j}.
\end{equation*}
Moreover the random variables $\Lambda^S(A)$, $A \subseteq V\backslash S$,
are independent. 

\subsection{Rule for Randomized Collapse}
Recall that the sequence of vertices chosen to collapse a hypergraph is
unimportant, provided we keep going until there are no more patches. 
However we shall use a specific randomized rule which turns out to
admit a description in terms of a finite-dimensional Markov chain.
This leads to a randomized process of collapsing hypergraphs $(\Lambda_n)_{n
\ge 0}$.  This will prove to be an effective means to compute the
numbers of identifiable vertices and identifiable hyperedges for $\Lambda_0$.

The process $(\Lambda_n)_{n \ge 0}$, together with a sequence of sets 
$(S_n)_{n \ge 0}$ such that $\Lambda_n=\Lambda^{S_n}$,
is constructed as follows. Let $S_0=\emptyset$ and $\Lambda_0=\Lambda$. 
Suppose that $S_n$ and $\Lambda_n$ have been defined.
If there are no patches in $\Lambda_n$, then $S_{n+1}=S_n$ and
$\Lambda_{n+1}=\Lambda_n$.
If there are patches in $\Lambda_n$, select one uniformly at random and denote 
by $v_{n+1}$ the corresponding vertex; then set $S_{n+1}=S_n\cup\{v_{n+1}\}$
and $\Lambda_{n+1}=\Lambda^{S_{n+1}}$.

\subsection{An embedded Markov chain}
Let $Y_n$ denote the number of patches and $Z_n$ the amount of debris in $\Lambda_n$.
Then $Y_n=0$ and $Z_n=|\Lambda^*|$ for $n \geq |V^*|$.
Also $|V^*| = \inf \{n \geq 0 \ : \ Y_n = 0\}$.
Let $W_{n+1}$ denote the number of other patches at time $n$ sharing the same 
vertex as 
the $(n+1)$st selected patch, and let $U_{n+1}$ denote the number of $2$-edges 
at time $n$ containing the $(n+1)$st selected vertex $v_{n+1}$.
Our analysis will rest on the observation that $(Y_n,Z_n)_{n \geq 0}$ is a
Markov chain, where, conditional on $Y_n = m \geq 1$ and $Z_n=k$, we have
$$
Y_{n+1} = Y_n-1-W_{n+1} + U_{n+1},\qquad Z_{n+1} = Z_n + 1 +W_{n+1}
$$
and where $W_{n+1} \sim B\left(m-1, 1/(N-n)\right)$ and $U_{n+1} \sim P\left((N-n-1)
\lambda_2 (N,n)\right)$ with $W_{n+1}$ and $U_{n+1}$ independent.

To see this, introduce the filtration
$$
{\cal F}_n = \sigma (S_r, Y_r, Z_r :  r=0,1,\dots, n).
$$
\begin{lemma}\label{l41}
Let 
$$
p(\lambda |S,k,m)=\PP[\Lambda^{S} = \lambda | \sum_{v} \Lambda^{S}(\{v\}) = m, 
\Lambda^{S}(\emptyset) = k].
$$
Then
$$
\PP[\Lambda_{n} = \lambda | {\cal F}_n] = p(\lambda |S_{n}, Y_{n}, Z_{n}).
$$
Equivalently, for all $B\in {\cal F}_n$ so that $B \subset \{S_n=S, Y_n = m, Z_n=k\}$,
$$
\PP[\Lambda_{n} = \lambda, B] = p(\lambda |S, m, k) \PP[B].
$$
\end{lemma}
The claimed Markov structure for $(Y_n,Z_n)_{n \geq 0}$ follows easily.


\begin{proof}
The identity is obvious for $n=0$.  Suppose it holds for $n$.
Let $B \subset \{S_n=S, Y_n = m, Z_n=k\}$.
Take $x\in V\backslash S$, $m'\geq 1$
and $k'> k$.  Set $S'=S\cup \{x\}$ and $B'=\{S_{n+1}=S', Y_{n+1} = m',
Z_{n+1}=k'\}\cap B$.  It will suffice to show, for all hypergraphs $\lambda'$
having $m'$ patches and amount of debris $k'$, that
$$
\PP[\Lambda_{n+1} = \lambda', B'] \propto 
p(\lambda' |S', m', k'),
$$
where $\propto$ denotes equality up to a constant independent of $\lambda'$.
But
$$
\PP[\Lambda_{n+1} = \lambda', B']  = \sum_{\lambda} \frac{k'-k}{m}
\PP[\Lambda_n=\lambda, B]
\propto \sum_\lambda  p(\lambda |S, m, k).
$$
where the sum is over all hypergraphs $\lambda$ which collapse to
$\lambda'$ on removing the vertex $x$. 
Let $Y^S, Z^S$ denote the number of patches, amount of debris in
$\Lambda^S$ respectively.
Since $Y^{S}$ and $Z^{S}$ are conditionally independent of 
$\Lambda^{S'}$ given $Y^{S'}$ and $Z^{S'}$, 
$$
\sum_\lambda  p(\lambda |S, m, k)
=\PP(\Lambda^{S'}=\lambda'|Y^S=m,Z^S=k)
\propto  p(\lambda' |S', m', k'),
$$
as desired.
\end{proof}

\section{Exponential Martingales for Jump Processes}
We recall here some standard notions for pure jump Markov processes in
$\mathbb{R}^d$ and their associated martingales.  These will be used to study
the fluid limit of a sequence of such jump processes in Section \ref{s6}.
\subsection{Laplace Transforms}
Let $(X_t)_{t\ge 0}$ be a pure jump Markov process taking values
in a subset $I$ of $\mathbb{R}^d$, with L\'evy kernel $K$.  Consider the
Laplace transform
\begin{equation*}
m(x,\theta) = \int_{\mathbb{R}^d} e^{\langle \theta, y \rangle} K(x,dy),
\quad\theta\in(\mathbb{R}^d)^*
\end{equation*}
and assume that, for some $\eta_0>0$,
\begin{equation}\label{C}
\sup_{x\in I} \sup_{|\theta|\le \eta_{0}} m(x,\theta)  \le  C < \infty
\end{equation}
The distribution of the time $T$ and displacement $\Delta X_T$ of the
first jump of $(X_t)_{t\ge 0}$ is given by
\begin{equation*}
\mathbb{P} (T\in dt, \ \Delta X_T \in dy \mid T>t, \ X_0=x)=K(x,dy)dt .
\end{equation*}
Introduce random measures $\mu$ and $\nu$ on $(0,\infty) \times \mathbb{R}^d$,
given by
\begin{align*}
\mu & = \sum_{\Delta X_t\not= 0} \ve_{(t,\Delta X_t)},\\
\nu (dt, dy) & = K(X_{t-}, dy) dt
\end{align*}
where $\ve_{(t,y)}$ denotes the unit mass at $(t,y)$; $\nu$ is thus the compensator of the random measure $\mu$, in the sense of \cite{K97}, p. 422.
\subsection{Martingales Associated with Jump Processes}
The fact that $\nu$ is a compensator implies that, 
for any previsible process
$a: \Omega \times(0,\infty) \times \mathbb{R}^d \to \mathbb{R}$ 
satisfying
\begin{equation*}
\mathbb{E} \int_{\mathbb{R}^{d}} |a(s,y)|\nu(ds, dy)<\infty,
\end{equation*}
the following process is a martingale
\begin{equation*}
\int^t_0 \int_{\mathbb{R}^{d}} a(s,y)(\mu-\nu) (ds,dy).
\end{equation*}
In particular, (\ref{C}) allows us to take $a(s,y)=y$, which gives the martingale
\begin{equation*}
M_t = \int^t_0 \int_{\mathbb{R}^{d}} y (\mu-\nu)(ds,dy).
\end{equation*}
Fix $\eta \in (0,\eta_0)$. Then there exists $A<\infty$ such that
\begin{equation}\label{A}
|m'' (x,\theta)|  \le A, \qquad x \in I, \quad |\theta|   \le
\eta, 
\end{equation}
where $'$ denotes differentation in $\theta$.  Define for $\theta \in (\mathbb{R}^d)^*$
\begin{equation*}
\phi(x,\theta) = \int_{\mathbb{R}^d} \{e^{\langle \theta,y \rangle} - 1 - \langle \theta, y \rangle \}K(x,dy).
\end{equation*}
Then $\phi \ge 0$ and, for $|\theta|  \le  \eta$, by the
second-order mean value theorem,
\begin{equation*}
\phi (x,\theta) = \int^1_0 m'' (x,r\theta)(\theta, \theta) (1-r) dr
\end{equation*}
so
\begin{equation*}
\phi (x,\theta)  \le  \frac{1}{2} A |\theta|^2, \qquad x \in
I, \quad |\theta | \le  \eta .
\end{equation*}
Let $(\theta_t)_{t\ge 0}$ be a previsible process in $(\mathbb{R}^d)^*$
with $|\theta_t|  \le  \eta$ for all $t$.  Set
\begin{equation}\label{expmart}
Z_t = Z^\theta_t = \exp \{\int^t_0 \langle\theta_s, dM_s\rangle -
\int^t_0\phi (X_s, \theta_s)ds\}.
\end{equation}
Then $(Z_t)_{t \ge 0}$ is locally bounded, 
and by the Dol\'{e}ans formula (\cite{K97}, p. 440),
\begin{equation*}
Z_t = 1 + \int^t_0 \int_{\mathbb{R}^d} Z_{s-} (e^{\langle \theta_{s},y \rangle} -1)
(\mu-\nu)(ds, dy).
\end{equation*}
Hence $(Z_t)_{t\ge 0}$ is a non-negative local martingale, so $\mathbb{E}
(Z_t) \leq  1$ for all $t$.  Hence
\begin{align*}
\mathbb{E} &\int^t_0 \int_{\mathbb{R}^d}| Z_{s-}
(e^{\langle \theta_s,y \rangle} - 1 )| \nu(ds,dy)\\
&\le \mathbb{E} \int^t_0 Z_s( m(X_s, \theta_s)+ m(X_s, 0))ds \le 2Ct 
\end{align*}
so $(Z_t)_{t\ge 0}$ is a martingale.
\begin{prop}\label{p51}
For all $\delta \in (0,A\eta t\sqrt{d} ]$
\begin{equation*}
\mathbb{P} \left(\sup_{s  <  t} |M_s|>\delta\right)  \le  (2d) e^{-\delta^{2}/(2Adt)}.
\end{equation*}
\end{prop}
\begin{proof}
Fix $\theta \in (\mathbb{R}^d)^*$ with $|\theta|=1$ and consider the stopping time
\begin{equation*}
T=\inf \{ t \ge 0: \langle \theta, M_t\rangle >\delta\}.
\end{equation*}
For $\varepsilon  <  \eta$, taking $\theta_t = \theta$ for all $t$
above, we know that $(Z^{\varepsilon \theta}_t)_{t \ge 0}$ is a
martingale.  On the set $\{T  \le  t\}$ we have $Z^{\varepsilon
\theta}_T \ge e^{\delta \varepsilon - At\varepsilon^{2}/2}$.  By
optional stopping
\begin{equation*}
\mathbb{E} (Z^{\varepsilon \theta}_{T\wedge t}) = \mathbb{E} (Z^{\varepsilon
\theta}_0) = 1.
\end{equation*}
Hence, 
\begin{equation*}
\mathbb{P} \left(\sup_{s \le  t} \ \langle\theta, M_s\rangle >\delta\right)
= \mathbb{P} (T  \le  t)  \le  e^{-\delta \varepsilon + At\varepsilon^{2}/2}.
\end{equation*}
When $\delta  \le  At\eta$ we can take $\varepsilon = \delta/At$ to obtain
\begin{equation*}
\mathbb{P} (\sup_{s  \le  t} \langle \theta, M_s\rangle> \delta )  \le  e^{-\delta^{2}/2At}.
\end{equation*}
Finally, if $\sup_{s  \le  t} |M_s|>\delta$, then $\sup_{s  \le 
t}\langle \theta, M_s\rangle > \delta/\sqrt{d}$ for one of $\theta = \pm
e_1,\dots, \pm e_d$.
\end{proof}

\section{Fluid Limit for Stopped Processes}\label{s6}
In this section we develop some general criteria for the convergence
of a sequence of Markov chains in $\mathbb{R}^d$ to the solution of a
differential equation, paying particular attention to the case where
the chain may stop abruptly on leaving a given open set.
\subsection{Fluid Limits}
It is possible to give criteria for the convergence of Markov processes
in terms of the limiting behaviour of their infinitesimal characteristics.
This is a powerful technique which has been intensively studied by
probabilists. The book of Ethier and Kurtz \cite{EK} is a key reference.
Further results are given in Chapter 17 of \cite{K97} and 
in Theorems IX.4.21 and IX.4.26 of \cite{JS87}. A particular case with
many applications is where the limiting process is deterministic and is
given by a differential equation, sometimes called a fluid limit.
The relevant probabilistic literature, though well developed, may not
be readily accessible to non-specialists seeking to apply the results in 
other fields. One field where fluid limits of Markov processes are
beginning to find interesting applications is random combinatorics.
Wormald \cite{W95} and co-workers have put forward a set of criteria which
is specially adapted to this application. The material in this section
may be considered as an alternative framework, somewhat more rigid but, we
hope, easy to use, developed with the same applications in mind.

Let $(X^N_t)_{t\ge 0}$ be a sequence of pure jump Markov processes
in $\mathbb{R}^d$. It may be that $(X^N_t)_{t\ge 0}$ takes values in some
discrete subset $I^N$ of $\mathbb{R}^d$ and that its L\'evy kernel $K^N(x,dy)$
is given naturally only for $x\in I^N$. 
So let us suppose that $I^N$ is measurable, that $(X^N_t)_{t \geq 0}$
takes values in
$I^N$, and that the  L\'evy kernel $K^N(x,dy)$ is given for $x\in I^N$.
Let $S$ be an open set in $\mathbb{R}^d$ and set $S^N=I^N\cap S$. 
We shall study, under
certain hypotheses, the limiting behaviour of $(X^N_t)_{t\ge0}$ as
$N\to\infty$, on compact time intervals,
up to the first time the process leaves $S$. 
In applications, the set $S$ will be chosen as the intersection of
two open sets $H$ and $U$. Our sequence of processes may all stop abruptly
on leaving some open set $H$, so that $K^N(x,dy)=0$ for $x\notin H$. If
this sort of behaviour does not occur, we simply take $H=\mathbb{R}^d$.
We choose $U$ so that the conjectured fluid limit path does not leave $U$
in the relevant compact time interval. Subject to this restriction we are free
to take $U$ as small as we like to facilitate the checking of convergence
and regularity conditions, which are required only on $S$.

The scope of our
study is motivated by the particular model which occupies the remainder
of this paper: so we are willing to impose a relatively strong, large
deviations-type, hypothesis on the L\'evy kernels $K^N$, see (\ref{A1}) below,
and we are interested to find that strong conclusions may be drawn using
rather direct arguments. On the other hand, in certain cases of our
model, the fluid limit path grazes the boundary of the set $S$: this
calls for a refinement of the usual fluid limit results to determine the
limiting distribution of the exit time.
\subsection{Assumptions}
Consider the Laplace transform
\begin{equation*}
m^N(x,\theta) = \int_{\mathbb{R}^{d}} e^{\langle\theta,y\rangle}K^N(x,dy),\quad
x\in S^N,\theta \in(\mathbb{R}^d)^*.
\end{equation*} 
We assume that, for some $\eta_0 >0$,
\begin{equation}\label{A1}
\sup_N\sup_{x\in
S^{N}}\sup_{|\theta| \le \eta_{0}}\dfrac{m^N(x,N\theta)}{N}< \infty.
\end{equation}
Set $b^N(x) = m^{N\prime}(x,0)$, where $'$ denotes the derivative in $\theta$.
We assume that, for some Lipschitz vector field $b$ on $S$,
\begin{equation}\label{A2}
\sup_{x\in
S^{N}}\quad\left|b^{N}(x) - b(x)\right| \to 0.
\end{equation}
We write $\tilde{b}$ for some Lipschitz vector field on
$\mathbb{R}^d$ extending $b$. (Such a $\tilde{b}$ is given, for example, by
$\tilde{b}(x)=\sup\{b(y)-K|x-y|:y\in S\}$ where $K$ is the Lipschitz constant
for $b$.) Fix a point $x_0$ in the closure $\bar{S}$ of $S$ and denote 
by $(x_t)_{t\ge 0}$ the unique solution to $\dot{x}_t =\tilde{b}(x_t)$ 
starting from $x_0$.
We assume finally that, for all $\delta>0$,
\begin{equation}\label{A4}
\limsup_{N\to\infty} {N^{-1}}\log\mathbb{P}(| {X^N_0} -x_0 |>\delta) < 0.
\end{equation}
Whilst these are not the weakest conditions for the fluid limit, they are 
readily verified in many examples of interest. In particular we will be able to 
verify them for the Markov chains associated with hypergraph collapse in
Section 3.
\subsection{Exponential convergence to the fluid limit}
Fix $t_0 > 0$ and set
\begin{equation*}
T^N =\inf \{t \ge 0 : X^N_t\notin S\} \wedge t_0.
\end{equation*} 

\begin{prop}\label{p61}
Under assumptions (\ref{A1}), (\ref{A2}), (\ref{A4}),
 we have, for all $\delta > 0$,
\begin{equation}\label{C1}
\limsup_{N\to\infty}{N^{-1}}\log\mathbb{P}\big(\sup_{t\le T^N}\big|{X^N_t -
x_t}\big|>\delta\big)<0. 
\end{equation}
\end{prop}

\begin{proof}  The following argument is widely known but we have not
found a convenient reference.  Set $b^N(x) = m^{N\prime}(x,0)$ and
define $(M^N_t)_{t\ge 0}$ by
\begin{equation*}
X^N_t = X^N_0 + M^N_t+ \int^t_0 b^N(X^N_s)ds.
\end{equation*} 
Note that $M^N_t$ corresponds to the martingale we identified in Proposition \ref{p51}.
Fix $\eta \in (0,\eta_0)$.  
Assumption (\ref{A1}) implies that
there exists $C<\infty$ such that, for all $N$,
\begin{equation*}
\big|m^{N\prime\prime}(x,\theta)\big|  \le  C/N,\qquad x\in S^N,
\quad |\theta|  \le  N\eta.
\end{equation*}
Compare this estimate with (\ref{A}). By applying Proposition \ref{p51}
to the stopped process $(X^N_{t\wedge T^{N}})_{t\ge 0}$, we find
constants $\varepsilon_0>0$ and $C_0<\infty$, depending only on
$C,\eta,d$ and $t_0$ such that, for all $N$ and all
$\varepsilon \in  (0,\varepsilon_0]$,
\begin{equation}\label{C2}
\mathbb{P} (\sup_{t  \le  T^{N}} |M^N_t| >\varepsilon)  \le  C_0
e^{-N\varepsilon^{2}/C_{0}}.
\end{equation}
Given $\delta>0$, set $\varepsilon = \min \{e^{-Kt_{0}}\delta/3,
\varepsilon_0\}$, where $K$ is the Lipschitz constant of $\tilde{b}$.  Let
\begin{equation*}
\Omega^N = \{|X^N_0 - x_0 |   \le  \varepsilon \ \hbox{and} \
\sup_{t  \le  T^{N}} |M^N_t |  \le  \varepsilon \}.
\end{equation*}
Then (\ref{A4}) and (\ref{C2}) together imply that
\begin{equation*}
\limsup_{N\to\infty} N^{-1} \log \mathbb{P} (\Omega\backslash \Omega^N)<0.
\end{equation*}
On the other hand, by (\ref{A2}), there exists $N_0$ such that
$|b^N(x)-b(x)|  \le \varepsilon/t_0$ for all $x\in S^N$ and all $N
\ge N_0$.  We note that
\begin{align*}
X^N_t - x_t & = (X^N_0 - x_0) + M^N_t + \int^t_0 (b^N (X^N_s)-b(X^N_s))ds\\
&+ \int^t_0 (\tilde b(X^N_s)-\tilde b(x_s))ds
\end{align*}
so, for $N\ge N_0$, on $\Omega^N$, for $t  \le  T^N$,
\begin{equation*}
|X^N_t-x_t|  \le  3\varepsilon + K \int^t_0 |X^N_s - x_s|ds
\end{equation*}
which implies, by Gronwall's lemma, that
$\sup_{t  \le  T^{N}} |X^N_t - x_t|  \le  \delta$.
\end{proof}
\subsection{Limiting Distribution of the Exit Time}
The remainder of this section is concerned with the question, left open
by Proposition \ref{p61}, of determining the limiting distribution of $T^N$.  Set
\begin{align*}
\tau & =  \inf \{ t \ge 0 :  x_t \notin \bar{S}\} \wedge t_0, \\
\mathcal{T} & = \{t\in [0,\tau) : x_t \notin S\}.
\end{align*}
It is straightforward to deduce from (\ref{C1}) that, for all $\delta >0$,
\begin{equation}\label{C3}
\limsup_{N\to\infty} N^{-1} \log \mathbb{P} \left(\inf_{t\in \mathcal{T} \cup\{\tau
\}} |T^N - t|>\delta\right) < 0 .
\end{equation}
In particular, if $\mathcal{T}$ is empty, then $T^N\to \tau$ in probability
and, for all $\delta>0$,
$$
\limsup_{N\to\infty}{N^{-1}}\log\mathbb{P}\big(\sup_{t\le t_0}\big|{X^N_t -
x_{t\wedge\tau}}\big|>\delta\big)<0.
$$ 
The reader who wishes only to know the proof of Theorem \ref{t1}
may skip to Section \ref{s7} as the remaining results of this section are
needed only for the more general case considered in Theorem \ref{t2}.
\subsection{Fluctuations}
We assume here that
\begin{equation}\label{A5}
\mathcal{T} \ \hbox{is finite} .
\end{equation}
In this case the limiting distribution of $T^N$ may be obtained from
that of the fluctuations $\gamma^N_t = \sqrt{N} (X^N_{t\wedge T^{N}} -
x_{t\wedge T^{N}})$. 
We assume that there exists a limit kernel $K(x,dy)$, defined for $x\in S$
such that $m(x,\theta)<\infty$ for all $x\in S$ and $|\theta| \le \eta_0$, 
where
\begin{equation*}
m(x,\theta) = \int_{\mathbb{R}^{d}} e^{\langle\theta,y\rangle}K(x,dy),\quad
x\in S, \quad\theta\in (\mathbb{R}^d)^*.
\end{equation*}
For convergence of the fluctuations we assume
\begin{equation}\label{A6}
\gamma^N_0 \to \gamma_0 \ \hbox{in distribution},
\end{equation}
\begin{equation}\label{A6a}
\sup_{x\in S^{N}}
\sup_{|\theta| \le \eta_{0}}\quad\left|\dfrac{m^N(x,N\theta)}{N}
- m(x,\theta)\right|\to 0,
\end{equation}
\begin{equation}\label{A7}
\sup_{x\in S^{N}} \sqrt{N} |b^N (x) - b(x)| \to 0,
\end{equation}
\begin{equation}\label{A8}
a \ \text{is Lipschitz and} \ b \ \text{is} \ C^1 \ \text{on} \ S,
\end{equation}
where $b^N(x) = m^{N\prime}(x,0)$ and $a(x) = m''(x,0)$. 
Of course (\ref{A6a}) will force $b(x) = m^\prime(x,0)$.
\subsection{Limiting Stochastic Differential Equation}
Consider the process $(\gamma_t)_{t  \le {\tau}}$ given by the
linear stochastic differential equation
\begin{equation}\label{E1}
d\gamma_t = \sigma(x_t) dB_t + \nabla b(x_t) \gamma_t dt
\end{equation}
and starting from $\gamma_0$, where $B$ is a Brownian motion and
$\sigma(x)\sigma(x)^*=a(x)$.  The distribution of $(\gamma_t)_{t
 \le  {\tau}}$ does not depend on the choice of $\sigma$.  For
convergence of $T^N$ we assume, in addition, 
\begin{equation}\label{A9}
\partial S \hbox{ is } C^1 \hbox{ at } x_t \hbox{ with inward normal } n_t,
\hbox{ and } \mathbb{P} (\langle n_t, \gamma_t\rangle = 0)=0,\hbox{ for all }
t\in \mathcal{T}.
\end{equation}
\begin{thm}\label{t3}
Under assumptions (\ref{A4}), (\ref{A5}), (\ref{A6}), (\ref{A6a}), (\ref{A7}), (\ref{A8}), (\ref{A9}) 
we have $T^N\to T$ in distribution, where
\begin{equation*}
T=\min \{t\in \mathcal{T} : \langle n_t, \gamma_t\rangle < 0\}\wedge \tau .
\end{equation*}
\end{thm}

\begin{proof}
Let $\tau_0=0$ and write the positive elements of $\mathcal{T}$ as
$\tau_1<\dots < \tau_m$.  
Define, for $k=0,1,\dots ,m$,
\begin{equation*}
\tilde{\gamma}^N_k  = 
\begin{cases} \gamma^N_{\tau_{k}} & \text{if} \ T^N>\tau_k\\
\partial & \text{otherwise} ,
\end{cases}
\end{equation*}
\begin{equation*}
\tilde{\gamma}_k = \begin{cases} \gamma_{\tau_{k}} & \text{if} \ T>\tau_k\\
\partial & \text{otherwise}, 
\end{cases}
\end{equation*}
where $\partial$ is some cemetery state.  We will show by induction, for
$k=0,1,\dots, m$, that
\begin{equation}\label{C4}
(\tilde{\gamma}^N_0,\dots , \tilde{\gamma}^N_k)\to
(\tilde{\gamma}_0,\dots , \tilde{\gamma}_k) \ \text{in
distribution}.
\end{equation}
Given (\ref{C3}), this implies that $T^N\to T$ in distribution, as required.

Note that both $(\tilde{\gamma}^N_k)_{0  \le  k \le   m}$ and
$(\tilde{\gamma}_k)_{0  \le  k  \le  m}$ may be considered as
time-dependent Markov processes.  Hence, by a conditioning argument, it
suffices to deal with the case where $\gamma_0$ is non-random. By (\ref{A9}),
if $x_0 \in \partial S$, we can assume that $\partial S$ is $C^1$ at
$x_0$ and $\langle n_0,\gamma_0\rangle \neq 0$.  Moreover, for the
inductive step, it suffices to consider the case where
$\tilde{\gamma}_k$ is non-random, not $\partial$, and to show that, if
$\tilde{\gamma}^N_k\to \tilde{\gamma}_k$ in probability, then
$\tilde{\gamma}^N_{k+1}\to \tilde{\gamma}_{k+1}$ in distribution.  We
lose no generality in considering only the case $k=0$.

We have assumed that $\gamma^N_0 \to \gamma_0$ in distribution.  Note
that $T=0$ if and only if $x_0 \in \partial S$ and $\langle n_0,
\gamma_0 \rangle < 0$.  On the other hand, since $X^N_0 = x_0 +
\sqrt{N}\gamma^N_0$, we have $\mathbb{P}(T^N=0)\to 1$ if and only if $x_0 \in
\partial S$ and $\langle n_0, \gamma_0 \rangle <0$.  Hence 
$\tilde{\gamma}^N_0\to \tilde{\gamma}_0$ in distribution, that is, (\ref{C4})
holds for $k=0$.

In Lemmas \ref{l65}, \ref{l66} and \ref{l67} below, we will show that, if $x_0 \in S$, or
$x_0 \in \partial S$ and $\langle n_0, \gamma_0\rangle >0$, then 
\begin{equation*}
\mathbb{P} (T^N>\varepsilon)\to 1 \hbox{ for some } \varepsilon >0, 
\end{equation*}
and, in the case $m \ge 1$,
\begin{align*}
& \gamma^N_{\tau_{1}}  \to \gamma_{\tau_{1}} \ \text{in distribution}, \\
&\mathbb{P} (\langle n_{\tau_{1}}, \gamma^N_{\tau_{1}} \rangle \ge 0 \
\text{and} \ T^N  \le  \tau_1)\to 0,\\
& \mathbb{P} (\langle n_{\tau_{1}}, \gamma^N_{\tau_{1}} \rangle < 0 \
\text{and} \ T^N > \tau_1) \to 0.
\end{align*}
It follows that $\tilde{\gamma}^N_1\to \tilde{\gamma}_1$ in
distribution, so (\ref{C4}) holds for $k=1$.  This establishes the induction
and completes the proof.
\end{proof}

We remark that the same proof applies when the L\'evy kernels $K^N$ have
a measurable dependence on the time parameter $t$, subject to obvious
modifications and to each hypothesis holding uniformly in $t \le  t_0$.

For the remainder of this section, the assumptions of Theorem \ref{t3} are in
force and $\gamma_0$ is non-random.

\begin{lemma}\label{l63}
For all $\varepsilon >0$ there exists $\lambda <\infty$ such that, for
all $N$
\begin{equation*}
\mathbb{P} (\sup_{t  \le  t_{0}} |\gamma^N_t | \ge \lambda) < \varepsilon.
\end{equation*}
\end{lemma}

\begin{proof} Given $\varepsilon>0$, choose $\lambda < \infty$ and $N_0$
such that, for
$\lambda' = e^{-Kt_{0}} \lambda/3$ and $N \ge N_0$
\begin{equation*} 
\sqrt{N} |b^N(x) - b(x)|  \le  \lambda'/t_0, \quad x \in S^N
\end{equation*}
and, with probability exceeding $1-\varepsilon$,
\begin{equation*}
|\gamma^N_0 |  \le  \lambda',
\end{equation*}
\begin{equation*}
\sqrt{N} \sup_{t  \le  T^{N}} |M^N_t|  \le  \lambda' .
\end{equation*}
This is possible by (\ref{C2}) and (\ref{A7}).  These three inequalities imply
\begin{equation*}
|\gamma^N_t|  \le  3\lambda' + K \int^t_0 |\gamma^N_s | ds, \qquad t
 \le  T^N,
\end{equation*}
so, by Gronwall's lemma
\begin{equation*}
\sup_{t  \le  T^{N}} |\gamma^N_t|  \le  \lambda.
\qed
\end{equation*}
\end{proof}

\begin{lemma}\label{l64}
For all $\varepsilon>0$ there exists $\lambda < \infty$ such that, for
all $\delta >0$, there exists $N_\delta <\infty$ such that, for all $N
\ge N_\delta$ and all $t  \le  t_0$,
\begin{equation*}
\mathbb{P} \left(\sup_{s  \le  t_{0}, t  \le  s  \le  t + \delta}
|\gamma^N_s - \gamma^N_t|> \lambda \sqrt{\delta}\right) < \varepsilon
\end{equation*}
\end{lemma}

\begin{proof}
Consider first the case $t=0$.  Given $\varepsilon >0$, choose $\lambda
<\infty$ such that, for all $\delta >0$, there exists $N_\delta$ such
that, for $\lambda' = e^{Kt_0}\lambda/3$ and $N \ge N_\delta$
\begin{equation*}
\sqrt{N} |b^N(x) - b(x)|  \le  \lambda'/\sqrt{t_{0}},\qquad x\in S^N,
\end{equation*}
and, with probability exceeding $1-\varepsilon$,
\begin{equation*}
|\gamma^N_0 |  \le  \lambda' /K\sqrt{t_{0}},
\end{equation*}
\begin{equation*}
\sqrt{N} \sup_{t  \le  T^{N}\wedge \delta} |M^N_t|  \le  \lambda'\sqrt{\delta}.
\end{equation*}
This is possible by (\ref{C2}) and (\ref{A7}).  These three inequalities imply
\begin{equation*}
|\gamma^N_t - \gamma^N_0 |  \le  3 \lambda' \sqrt{\delta} + K
\int^t_0 |\gamma^N_s - \gamma^N_0 |ds, \qquad t  \le  T^N \wedge \delta,
\end{equation*}
so by Gronwall's lemma
\begin{equation*}
\sup_{t  \le  T^{N}\wedge \delta} |\gamma^N_t - \gamma^N_0 |
 \le  \lambda \sqrt{\delta}.
\end{equation*}
The case $t>0$ follows by the same sort of argument, using Lemma \ref{l63} to
get the necessary tightness of $\gamma^N_t$.
\end{proof}

\begin{lemma}\label{l65}
Suppose either $x_0\in S$, or $x_0 \in \partial S$ and $\langle n_0,
\gamma_0\rangle >0$.  Then $\mathbb{P} (T^N>\varepsilon)\to 1$ as $N\to \infty$
for some $\varepsilon>0$.
\end{lemma}

\begin{proof}
The case $x_0 \in S$ follows from (\ref{C3}).  Suppose then that $x_0\in
\partial S$ and $\langle n_0,\gamma_0 \rangle > 0$.  Then, since
$\partial S$ is $C^1$ at $x_0$, for all $\varepsilon >0$, there exists 
$\delta(\varepsilon) >0$ such that, for all $x\in \bar{S}$ with 
$|x-x_{0}| \le  \delta(\varepsilon)$, and all $v \in \mathbb{R}^d$,
\begin{equation}\label{G}
|v|  \le  \delta(\varepsilon) \ \text{and} \ \langle n_0, v\rangle
\ge \varepsilon |v| \Rightarrow \ x + v \in S.
\end{equation}
Since $\langle n_0,\gamma_0\rangle > 0$, by Lemma \ref{l64}, given
$\varepsilon >0$ there exist $\varepsilon_1 >0$ and $N_0$ such that,
for all $N \ge N_0$ and $t  \le  T^N \wedge \varepsilon_1$,
\begin{equation*}
\langle n_0, \gamma^N_t\rangle >\varepsilon_1 |\gamma^N_t|,\quad
|\gamma^N_t|< 1/\varepsilon_1,
\end{equation*}
with probability exceeding $1-\varepsilon$.  Choose $\varepsilon_2 \in
(0,\varepsilon_1)$ so that $|x_t-x_0|  \le   \delta (\varepsilon_1)$
and $x_t \in \bar{S}$ whenever $t  \le  \varepsilon_2$. Set $N_1 =
\max \{N_0, (\varepsilon_1\delta (\varepsilon_1))^{-2}\}$, then, 
for $N \ge N_1$ and $t  \le  T^N\wedge \varepsilon_2$,
\begin{equation}\label{C5}
x_t\in \overline{S}, \ |x_t-x_0| \le \delta(\varepsilon_1),\
N^{-1/2}|\gamma^N_t| \le \delta(\varepsilon_1), \ \langle 
n_0,\gamma^N_t\rangle>\varepsilon_1|\gamma^N_t|,
\end{equation} 
with probability exceeding $1-\varepsilon$. By (\ref{G}), (\ref{C5}) implies $X^N_t=
x_t + N^{-1/2}\gamma^{N}_t \in S$.   Hence $\mathbb{P}(T^N  \le
\varepsilon_2)<\varepsilon$ for all $N \ge N_1$.
\end{proof}
For the rest of this section we assume that $m \ge 1$.  (The next
result holds with $\tau_1$ replaced $\tau$ when $m=0$, by the same
argument, but we do not need this.)

\begin{lemma}\label{l66}
Suppose either $x_0 \in S$, or $x_0 \in \partial S$ and $\langle n_0,
\gamma_0 \rangle > 0$.  Then $\gamma^N_{\tau_{1}} \to \gamma_{\tau_{1}}$
in distribution as $N\to \infty$.
\end{lemma}

\begin{proof}
By Lemma \ref{l64}, given $\delta >0$, we can find $t<\tau_1$ such that, for
all $N$,
\begin{equation*}
\mathbb{P} (|\gamma^N_t - \gamma^N_{\tau_{1}} |>\delta )<\delta,\qquad  
\mathbb{P}(|\gamma_t - \gamma_{\tau_{1}} |>\delta)< \delta.
\end{equation*}
Hence it suffices to show $\gamma^N_t\to \gamma_t$ in distribution for
all $t<\tau_1$.

Define $(\psi_t)_{t  \le  \tau}$ in $\mathbb{R}^d \otimes (\mathbb{R}^d)^*$ by
\begin{equation*}
\dot\psi_t = \nabla b(x_t) \psi_t, \qquad \psi_0 = id.
\end{equation*}
Fix $\theta \in (\mathbb{R}^d)^*$ and set $\theta_t = (\psi^*_t)^{-1}\theta$. Then
\begin{equation*}
d\langle \theta_t, \gamma_t\rangle =  \langle \theta_t, \sigma (x_t)
dB_t\rangle, \qquad t  \le  \tau,
\end{equation*}
so
\begin{equation*}
\langle \theta_t, \gamma_t\rangle \sim N(\langle \theta,
\gamma_0\rangle, \int^t_0 \langle \theta_s, a(x_s)\theta_s \rangle
ds),\quad t  \le  \tau.
\end{equation*}
On the other hand, for $(M^N_t)_{t \ge 0}$ as in the proof of
Proposition \ref{p61},
\begin{equation*}
d\langle \theta_t, \gamma^N_t\rangle = \sqrt{N}\langle \theta_t,
dM^N_t\rangle + R^{N,\theta}_t dt, \qquad t  \le  T^N
\end{equation*}
where
\begin{equation*}
R^{N,\theta}_t = \sqrt{N} \langle \theta_t, b^N (X^N_t)- b(x_t) - \nabla
b(x_t)(X^N_t - x_t)\rangle.
\end{equation*}
By (\ref{A7}),
\begin{equation*}
\sup_{t  \le  T^{N}} \sqrt{N} |b^N (X^N_t)- b(X^N_t) |\to 0.
\end{equation*}
By (\ref{A8}), given $\varepsilon>0$,  there exists $\delta >0$ such that, for
all $t\in [\varepsilon, \tau_1-\varepsilon]$, for $|x-x_t|  \le  \delta$,
\begin{equation*}
|b(x)-b(x_t)-\nabla b(x_t)(x-x_t) |  \le  \varepsilon |x-x_t|.
\end{equation*}
Hence $|X^N_t - x_t|  \le  \delta$ and
 $\varepsilon  \le t \le  \tau_1 - \varepsilon$ imply
\begin{equation*}
\sqrt{N} |b(X^N_t) - b(x_t) - \nabla b(x_t) (X^N_t - x_t)|  \le  \varepsilon|\gamma^N_t|.
\end{equation*}
Combining this with Lemma \ref{l63}, we deduce that
\begin{equation*}
\int^{\tau_{1}}_0 |R^{N,\theta}_t |dt \to 0\quad \text{ in probability}.
\end{equation*}
Hence it suffices to show, for all $\theta \in (\mathbb{R}^d)^*$ and all $t<\tau_1$,
\begin{equation*}
\sqrt{N} \int^t_0 \langle \theta_s, dM^N_s\rangle \rightarrow N(0, \int^t_0
\langle \theta_s, a(x_s)\theta_s\rangle ds)\quad\hbox{in distribution.}
\end{equation*} 
Indeed, it suffices to show, for all $\theta \in (\mathbb{R}^d)^*$ and
$t<\tau_1$, that $\mathbb{E} (E^{N,\theta}_t)\to 1$ as $N\to \infty$, where
\begin{equation*}
E^{N,\theta}_t = \exp \lbrace i\sqrt{N}\int^t_0 \langle \theta_s, dM^N_s
\rangle + \frac{1}{2} \int^t_0 \langle \theta_s, a(x_s)\theta_s\rangle ds\rbrace.
\end{equation*}
Set $\tilde{m}^N(x,\theta) = m^N (x,i\theta)$,
$\tilde{m}(x,\theta)=m(x,i\theta)$ and
\begin{equation*}
\tilde{\phi}^N (x,\theta) = 
\int_{\mathbb{R}^{d}} (e^{i\langle \theta,y\rangle} - 1-i\langle \theta, y\rangle ) K^N (x,dy).
\end{equation*}
By (\ref{A6a}), for all $\eta<\eta_0$, we have
\begin{equation*}
\sup_{x\in S^{N}} \sup_{|\theta|  \le  \eta}
|N\tilde{m}^{N\prime\prime} (x,N\theta) - \tilde{m}'' (x,\theta)| \to 0.
\end{equation*}
Note that
\begin{align*}
&\tilde{\phi}^N (x,\sqrt{N}\theta) + \frac{1}{2} \langle \theta,
a(x)\theta\rangle \\
&=\int^1_0 \left(N\tilde{m}^{N\prime\prime}(x,\sqrt{N} r\theta) -
\tilde{m}'' (x,0)\right) (\theta, \theta)(1-r) dr 
\end{align*}
so, for all $\rho <\infty$,
\begin{equation}\label{C6}
\sup_{x\in S^N} \sup_{|\theta|  \le  \rho} |\tilde{\phi}^N (x,
\sqrt{N}\theta) + \frac{1}{2} \langle \theta, a(x) \theta \rangle | \to 0.
\end{equation}
Write $E^{N,\theta}_t = E^N_t = Z^N_t A^N_t B^N_t$, where
\begin{align*}
Z^N_t & = \exp \{i \sqrt{N} \int^t_0 \langle \theta_s, dM^N_s\rangle -
\int^t_0 \tilde{\phi}^N (X^N_s, \sqrt{N} \theta_s)ds\},\\
A^N_t & = \exp \{ \int^t_0 (\tilde{\phi}^N (X^N_s, \sqrt{N}\theta_s) +
\frac{1}{2} \langle \theta_s, a (X^N_s)\theta_s\rangle )ds \},\\
B^N_t & = \exp \{\int^t_0 \frac{1}{2} \langle \theta_s, (a (x_s) - a
(X^N_s))\theta_s \rangle ds\}.
\end{align*}

Now $(Z^N_{t\wedge T^{N}})_{t  \le{\tau}}$ is a martingale, as
in (\ref{expmart}), so $\mathbb{E}(Z^N_{t\wedge T^N})=1$ for all $N$.  
Fix $t  \le {\tau}$.
 By (\ref{C6}), $Z^N_{t\wedge T^{N}}$ is bounded, uniformly in
$N$, and $A^N_{t\wedge T^{N}}\to 1$ uniformly as $N\to \infty$. 
Moreover, by (\ref{A8}), $B^N_{t\wedge T^{N}}$ is bounded uniformly in $N$ and
converges to 1 in probability, using (\ref{C1}). Hence
\begin{equation*}
\mathbb{E}(Z^N_{t\wedge T^{N}} A^N_{t\wedge T^{N}} B^N_{t\wedge T^{N}})\to 1
\end{equation*}
as $N\to \infty$.  By Lemma \ref{l65} and (\ref{C3}), $\mathbb{P}(T^N>t)\to 1$ for all
$t<\tau_1$. It follows that $\mathbb{E}(E^N_t)\to 1$ for all $t<\tau_1$ as required.
\end{proof}

\begin{lemma}\label{l67}
Suppose either $x_0 \in S$, or $x_0 \in \partial S$ and $\langle n_0,
\gamma_0 \rangle >0$.  Then, as $N\to \infty$,
\begin{align*}
& \mathbb{P}(\langle n_{\tau_{1}}, \gamma^N_{\tau_{1}} \rangle \ge 0 
\text{ and } T^N  \le  \tau_1) \to 0, \\
& \mathbb{P}(\langle n_{\tau_{1}}, \gamma^N_{\tau_{1}} \rangle < 0  \text{ and
} T^N > \tau_1)\to 0.
\end{align*}
\end{lemma}

\begin{proof}
By Lemma \ref{l66}, given $\varepsilon>0$, there exists $\varepsilon_1>0$ and
$N_0$ such that, for all $N \ge N_0$
\begin{equation*}
|\langle n_{\tau_{1}}, \gamma^N_{\tau_{1}} \rangle |>\varepsilon_{1}
|\gamma^N_{\tau_{1}}|,\quad \varepsilon_{1}< |\gamma^N_{\tau_{1}} | < 1/\varepsilon_{1},
\end{equation*}
with probability exceeding  $1-\varepsilon$.  Then by Lemma \ref{l64}, there
exists $\varepsilon_2>0$ and $N_1 \ge N_0$ such that, for all $N
\ge N_1$, with probability exceeding $1-\varepsilon$, {\em
either\/} 
\begin{equation}\label{D1}
\langle n_{\tau_{1}}, \gamma^N_t \rangle >\varepsilon_2
|\gamma^N_t|,\quad |\gamma^N_t |<1/\varepsilon_2 \quad \text{ for all }
t \in [\tau_1-\varepsilon_2, \tau_1]
\end{equation}
{\em or}
\begin{equation}\label{D2}
\langle n_{\tau_{1}}, \gamma^N_{\tau_{1}} \rangle <-\varepsilon_{2}
|\gamma^N_{\tau_{1}}|,\quad |\gamma^N_{\tau_{1}} |<1/\varepsilon_{2}
.
\end{equation}
Since $\partial S$ is $C^1$ at $x_{\tau_{1}}$, there exists $\delta>0$
such that 
\begin{equation*}
\hbox{ if }x\in \bar{S}\hbox{ and }v\in \mathbb{R}^d\hbox{ with
}|x-x_{\tau_{1}}|  \le  \delta, |v|  \le  \delta\hbox{ and
}\langle n_{\tau_{1}},v\rangle < \varepsilon_2 |v|,\hbox{ then }x+v\in S 
\end{equation*}
and 
\begin{equation*}
\hbox {if }v\in \mathbb{R}^d\hbox{ with }|v|  \le  \delta\hbox{ and
}\langle n_{\tau_{1}}, v\rangle < - \varepsilon_2 |v|,\hbox{ then
}x_{\tau_{1}}+v \notin S. 
\end{equation*}
Choose $\varepsilon_3 \in (0,\varepsilon_2]$ such that
$|x_t-x_{\tau_{1}}|  \le  \delta$ and $x_t \in \bar{S}$ whenever
$t\in [\tau_1-\varepsilon_3, \tau_1]$.  Set $N_2 = \max \{N_1,
(\varepsilon_2\delta)^{-2}\}$. Then, for $N \ge N_2$, since $X^N_t
= x_t + N^{-1/2} \gamma^N_t$ on $\{T^N \ge t\}$, (\ref{D1}) implies
$X^N_t\in S$ for all $t\in [\tau_1 - \varepsilon_2, \tau_1]$ or
$T^N<\tau_1-\varepsilon_2$, and (\ref{D2}) implies $X^N_{\tau_{1}} \notin S$
or $T^N<\tau_1$.  We know by Lemma \ref{l65} and (\ref{C3}) that
$\mathbb{P}(T^N<\tau_1-\varepsilon_2)\to 0$ as $N\to \infty$.  Hence, with high
probability, as $N\to \infty$, $\langle n_{\tau_{1}},
\gamma^N_{\tau_{1}}\rangle \ge 0$ implies (\ref{D1}) and then
$T^N>\tau_1$, and $\langle n_{\tau_{1}}, \gamma^N_{\tau_{1}}\rangle <0$
implies (\ref{D2}) and then $T^N  \le  \tau_1$.
\end{proof}

\section{Fluid Limit of Collapsing Hypergraphs}\label{s7}
We now apply the general theory from the preceding sections to prove
our main results Theorems \ref{t1} and \ref{t2}.
\subsection{L\'evy Kernel for Collapse of Random Hypergraphs}
In Section \ref{s4} we introduced a Markov process $(\Lambda_n)_{n\ge 0}$ of
collapsing hypergraphs, starting from $\Lambda_0\sim$ Poisson$(\beta )$
and stopping when $n=|V^*|$, the number of identifiable vertices in
$\Lambda_0$. The process $(Y_n, Z_n)_{n\ge 0}$ of patches and debris in
$\Lambda_n$ was found itself to be Markov. We now view this process as a
function of the initial number of vertices $N$ and obtain a fluid limit
result when $N\to \infty$.

It will be convenient to embed our process in continuous time, by
removing vertices according to a Poisson process $(\nu_t)_{t\ge
0}$ of rate $N$ which stops when $\nu_t = | V^*|$. Set
\begin{equation*}
X^N_t = N^{-1}(\nu_t,Y_{\nu_{t}}, Z_{\nu_{t}})
\end{equation*} 
and note that $X^N$ takes values in
\begin{equation}\label{IN}
I^N = \{x\in \mathbb{R}^3 : Nx^1\in \{0,1,\dots,N-1\}, Nx^2, Nx^3 \in \mathbb{Z}^+\}
\cup\{(1,0,x^3) : Nx^3 \in \mathbb{Z}^+\}.
\end{equation}
The L\'evy kernel $K^N(x,dy)$ for $(X^N_t)_{t\ge 0}$ is naturally
defined for $x \in I^N$. If $x^2 = 0$ then $K^N(x,dy) = 0$. If $x^2 >
0$, then $N^{-1}K^N(x, \cdot)$ is a probability measure; by Lemma
\ref{l41}, it is the law
of the random variable $J^N/ N$, where
\begin{equation*}
J^N = (1,-1-W^N+U^N, 1+W^N),
\end{equation*}
\begin{equation*}
 W^N\sim B(Nx^2 -1, \quad 1/(N-Nx^1) ),\quad  U^N\sim P((N-Nx^1-1)
\lambda_2 (N,Nx^1))
\end{equation*}
with $W^N$ and $U^N$ independent.

Recall that $R$ denotes the radius of convergence of the power series
$\beta(t)$, given by (\ref{beta}). We assume, until further notice, that
$R>0$ and fix $t_0\in(0,R\wedge 1)$ and $\rho\in (t_0,R\wedge 1)$.

\begin{lemma}\label{l71}
There is a constant $C < \infty$ such that
\begin{equation*}
| N \lambda_2 (N,n) - \beta''(n/N) |  \le   C(\log N)^2/N 
\end{equation*}
for all $N \in \mathbb{N}$ and $n \in \{0,1,\dots, [N\rho]\}$.
\end{lemma}

\begin{proof} 
Recall that
\begin{equation*}
\lambda_2 (N,n)  = N \sum_{i=0}^n (i+1)(i+2) \beta_{i+2}\quad
\frac{n(n-1)\dots(n-i+1)}{N(N-1)\dots(N-i-1)} .
\end{equation*}
Set $M  = A \log N$ where $A = (\log (R/\rho))^{-1} < \infty$. Then, for
$n  \le  [N\rho]$,
\begin{align*}
&|N\lambda_2(N,n) - \beta''(n/N)|\\
&  \le  \sum_{i=1}^{M\wedge n} (i+1)(i+2)\beta_{i+2}\delta_i(N,n) +
2\sum_{i=M+1}^\infty (i+1) (i+2) \beta_{i+2} \rho_i (N,n)
\end{align*}
where
\begin{equation*}
\delta_i (N,n)  = \bigg|\frac{N^2}{(N-i)(N-i-1)}\frac{n}{N}
\bigg(\frac{n-1}{N-1} \bigg)\dots \bigg(\frac{n-i+1}{N-i+1}\bigg) -
\bigg(\frac{n}{N}\bigg)^i \bigg|
\end{equation*}
and
\begin{equation*}
\rho_i (N,n) = \frac{N^2}{(N-i)(N-i-1)} \bigg(\frac{n}{N}\bigg)^i
 \le   C\rho^i.
\end{equation*}
Note that, for $j=0,\dots,i-1$ and $i\leq M\wedge n$,
$$
\left|\frac{n-j}{N-j}-\frac{n}{N}\right|\leq A\log N/N
$$
so, making use of the inequality $|\prod a_j -\prod b_j|\leq \sum |a_j-b_j|$
for $0\leq a_j, b_j\leq 1$, we obtain
$$
\delta_i (N,n)  \leq  C (\log N)^2  \rho^i/N.
$$
Hence
\begin{equation*}
| N\lambda_2 (N,n) - \beta''(n/N) |\  \le  \ C (\log N)^2 \beta''
(\rho) / N + C (\rho/R)^M
\end{equation*}
and $(\rho/ R)^M = 1/N$.
\end{proof}
\subsection{Fluid Limit}
The main result of this section is to obtain the limiting behaviour of
$(X^N_t)_{t\ge 0}$ as $N\to\infty$, which we deduce from Proposition
\ref{p61} and Theorem \ref{t3}.
 We present first the calculations by which the limit was discovered.

Note that, as $N\to \infty$, for $x^1<R\wedge 1$, we have $W^N\to W$ and
$U^N\to U$ in distribution, where
\begin{equation*}
W\sim P(x^2/(1-x^1)), \quad U \sim P((1-x^1)\beta'' (x^1)).
\end{equation*}
Set $J=(1,-1-W+U, \ 1+W)$.  Note also that $X^N_0\to x_0 = (0,\beta_1,
\beta_0) \text{ and } \sqrt{N}(X^N_0-x_0)\to \gamma_0$ in distribution,
where $\gamma^1_0=0, \gamma^2_0 \sim N(0,\beta_1)$, $\gamma^3_0 \sim
N(0,\beta_0)$, with $\gamma^2_0$ and $\gamma^3_0$ independent.  Thus,
subject to certain technical conditions, to be checked later, at least
up to the first time that $X^{N,1}_t \ge R\wedge 1$  or $X^{N,2}_t
= 0$, the limit path is given by $\dot{x}_t = b(x_t)$, starting from 
$x_0$, where
\begin{equation*}
b(x) = \mathbb{E}(J)=\left(1,-1-\frac{x^2}{1-x^1} + (1-x^1)\beta'' (x^1), \ \frac{x^2}{1-x^1}.\right)
\end{equation*}

Fix $\rho'\in (0,\infty)$ and set
\begin{equation}\label{statespace}
S=\{(x^1,x^2, x^3) : |x^1| <\rho, \quad x^2 \in (0,\rho'), \quad x^3\in \mathbb{R}\}
\end{equation}
then $b$ is Lipschitz on $S$ and, for $\rho'$ sufficiently large, the
maximal solution on $[0,t_0]$ to $\dot{x}_t = b(x_t)$ in $\bar{S}$
starting from $x_0$ is given by $(x_t)_{t \le\tau}$, where
\begin{equation*}
x_t = (t, (1-t)(\beta'(t)+\log (1-t)),\quad \beta(t)-(1-t)\log (1-t))
\end{equation*}
and
\begin{equation*}
\tau = z^* \wedge t_0.
\end{equation*}
\subsection{Limiting Fluctuations}
Set $a(x) = \mathbb{E}(J\otimes J)$.  A convenient choice of $\sigma$ such that
$\sigma \sigma^*=a$ is $\sigma = (V_1, V_2, V_3)$, where
\begin{equation*}
V_1(x) = \sqrt{\frac{x^2}{1-x^1}} 
\left(\begin{matrix} \phantom{-}0 \\
\phantom{-}1 \\-1\end{matrix}\right), 
\quad V_2 (x) =
\sqrt{(1-x^1)\beta''(x^1)}
\left(\begin{matrix} 0 \\ 1 \\0\end{matrix}\right), 
\quad V_3 (x) = b(x).
\end{equation*} 
Note that $a$ is a Lipschitz and $b$ is $C^1$ on $S$.  The limiting
fluctuations are given by
\begin{equation*}
d\gamma_t = \sum_{i} V_i  (x_t) dB^i_t +\nabla b(x_t)\gamma_t dt, \quad t
 \le  {\tau}
\end{equation*}
starting from $\gamma_0$, where $B$ is a Brownian motion in $\mathbb{R}^3$
independent of $\gamma_0$.
Note that
\begin{equation*}
\mathcal{T} = \{ t \in [0,{\tau}) : x_t \notin S\} = \zeta \cap [0,t_0).
\end{equation*}
In cases where $\mathcal{T}$ is non-empty, the limiting behaviour of
$(X^N_t)_{t  \le  t_{0}}$ depends on the signs of the component of
the fluctuations normal to the boundary, that is, on $(\gamma^2_t : t
\in \mathcal{T})$.

Note that $\theta_t=b(x_t) B^3_t$ satisfies
\begin{equation*}
d\theta_t = V_3 (x_t) dB^3_t + \nabla b(x_t) \theta_t dt.
\end{equation*}
This is the part of the fluctuations which reflects our Poissonization
of the time-scale.  Since $b^2(x_t)=0$ for all $t\in \mathcal{T}$, it does
not affect $(\gamma^2_t : t \in \mathcal{T})$.  So consider $\gamma^*_t =
\gamma_t - \theta_t$. Then
\begin{equation*}
d\gamma^*_t = V_1 (x_t) dB^1_t + V_2 (x_t) dB^2_t + \nabla b(x_t)
\gamma^*_t dt.
\end{equation*}
Note that $V^1_1(x) = V^1_2 (x) = 0$ and $\nabla b^1(x)=0$, so
$(\gamma^*_t)^1=(\gamma^*_0)^1=0$ for all $t$.  Also $\partial
b^2/\partial x^2 = - 1/(1-x^1)$ and $\partial b^2/\partial x^3=0$. Also
$x^1_t=t$ and $x^2_t/(1-x^1_t)=\beta'(t)+\log (1-t)$.  The sign of
$(\gamma^*_t)^2$ is the same as that of $\alpha_t =
(\gamma^*_t)^2/(1-t)$.  We have
\begin{align*}
d\alpha_t & = d \gamma_t^{*2}/ (1-t) + \gamma^{*2}_t / (1-t)^2 dt \\
& = (V^2_1 (x_t) dB^1_t+V^2_2 (x_t) dB^2_t)/(1-t) 
\end{align*}
so we can write $\alpha_{t} = W(\sigma^{2}_{t})$, where $W$ is a Brownian
motion and
\begin{align*}
\sigma^2_t & = \beta_1 + \int^t_0 \frac{\beta' (s)+\log (1-s)+(1-s)\beta''(s)}{1-s}ds\\
& = \frac{\beta'(t) + \log (1-t)+t}{1-t}.
\end{align*}
We have shown that $(\text{sgn }(\gamma^2_t):t\in \mathcal{T})$ has the same
distribution as $(\text{sgn}(W_{t/(1-t)}):t\in \mathcal{T}).$ In particular
$\mathbb{P}(\gamma^2_t = 0)=0$ for all $t\in \mathcal{T}$.

Recall that $Z$ is defined by
\begin{equation*}
Z=\min \{z\in \zeta : W \left(z/1-z\right)<0\} \wedge z^*.
\end{equation*}
Set
\begin{equation*}
T^N = \inf \{t \ge  0 : X^{N,2}_t = 0\}
\end{equation*}
and put $Z(t_0) = Z\wedge t_0$, $T^N(t_0) = T^N \wedge t_0$.

\begin{thm}\label{t4}
For all $\delta > 0$ we have
\begin{equation*}
\limsup_{N \to\infty}N^{-1}\log \mathbb{P}\left(\sup_{t  \le  T^N(t_{0})}
|X^N_t - x_t| > \delta\right)\quad < 0.
\end{equation*}
Moreover, $T^N(t_0) \to Z(t_0)$ in distribution as $N \to \infty$.
\end{thm}

\begin{proof}
We defined $I^N$, the state-space of $(X^N_t)_{t \ge 0}$, in (\ref{IN}),
 and $S$ in (\ref{statespace}). Set $S^N = I^N \cap S$.
For $x \in S^N$ we have
\begin{equation*}\
m^N(x,\theta) =\int_{\mathbb{R}^3} e ^{\langle\theta,y\rangle} K^N(x,dy) = N \mathbb{E}\left(e^{\langle\theta,J^N\rangle/N}\right)
\end{equation*}
so
\begin{align*}
\frac{m^N(x,\theta)}{N} & =
\exp\bigg\lbrace\theta_1-\theta_2+\theta_3+B\left(Nx^2-1,
\frac{1}{N-Nx^1}, \frac{\theta_3-\theta_2}{N}\right) \\
&+ P\left((N-Nx^1 -1) \lambda_2(N,Nx^1), \frac{\theta_2}{N}\right)\bigg\rbrace.
\end{align*}
where, for $\theta\in\mathbb{R}$, we write $B(N,p,\theta) = N \log (1-p+p e^\theta)$ and $P(\lambda,\theta) = \lambda(e^\theta-1)$.
So, by Lemma 6.1,
\begin{equation*}
\sup_{x\in S^{N}}\sup_{|\theta |\le\eta_{0}}
\bigg|\frac{m^N(x,N\theta)}{N} - m(x,\theta)\bigg|\to 0
\end{equation*}
as $N\to\infty$, for all $\eta_0 > 0$, where
\begin{equation*}
m(x,\theta) = \mathbb{E}(e^{\langle\theta,J\rangle}) =
\exp\{\theta_1-\theta_2+\theta_3+P\left(\frac{x^2}{1-x^1},\theta_3-\theta_2\right)+P\left((1-x^1)\beta''(x^1),\theta_2 \right)\}.
\end{equation*}
Set
\begin{equation*}
b^N(x) = \int_{\mathbb{R}^3} y \ K^N(x,dy)= \mathbb{E}(J^N)
\end{equation*}
then, by Lemma 7.1,
\begin{equation*}
\sup_{x \in S^{N}}  \sqrt{N} \ |b^N (x) - b(x)| \to 0.
\end{equation*}
Recall that
\begin{equation*}
X_0^{N,1} = 0,\quad NX_0^{N,2} \sim P(N\beta_1), \quad NX_0^{N,3} \sim P(N\beta_0)
\end{equation*}
and $x_0 = (0, \beta_1, \beta_0)$. By standard exponential estimates,
for all $\delta > 0$
\begin{equation*}
\limsup_{N\to\infty} \ N^{-1} \log\mathbb{P}(|X^N_0 - x_0 |>\delta) <0.
\end{equation*}
We have now checked the validity of (\ref{A2}), (\ref{A4}), (\ref{A5}), (\ref{A6}), (\ref{A6a}),
(\ref{A7}), (\ref{A8}), (\ref{A9}) in this context, 
so Proposition \ref{p61} and Theorem \ref{t3}
 apply to give the desired conclusions.
\end{proof}
\begin{rem}
If $z^*<1$, then $z^*<R \wedge 1$, so by choosing $t_0 \in
(z^*, R\wedge 1)$ we get $Z(t_0) = Z$ and, as $N\to\infty$, with high
probability $T^N(t_0) = T^N$.
Hence, when $z^*<1$, Theorem \ref{t4} holds with $Z$ and $T^N$ replacing
$Z(t_0)$ and $T^N(t_0)$. In particular, Theorem \ref{t1} follows.
\end{rem}
\subsection{Proof of Theorem \ref{t2}}
\begin{proof}
Recall that
\begin{equation*}
X^N_{T{^N}} = \bigg(\frac{|V^{N*}|}{N},\ 0, \ \frac{|\Lambda^{N*}|}{N}\bigg).
\end{equation*}
Let $z\in\zeta\cup\{z^*\}$. If $z<1$, then $z<R\wedge 1$ so, by choosing
$t_0\grave\in (z, R\wedge 1)$ in Theorem \ref{t3}, we obtain 
\begin{equation}\label{Z}
\mathbb{P}\bigg(| X^N_{T{^N}} - x_z |\le\delta\bigg)\to\mathbb{P}(Z=z)
\end{equation}
for all sufficiently small $\delta>0$

It remains to deal with the case $z = z^* =1$. Note that $|V^{N*}|\le N$
and $|\Lambda^{N*}| \le |\Lambda^N|$. Now $|\Lambda^N|\sim
P(N\beta(1))$ so $|\Lambda^N|/N\to \beta (1)$ in probability as
$N\to\infty$. It therefore suffices to show, for all $\delta >0$ and
$\alpha < \beta (1)-\delta$,
\begin{equation*}
\liminf_{N\to\infty}\quad\mathbb{P}\bigg(\frac{|V^{N*}|}{N}\ge
1-\delta\text{ and }\frac{|\Lambda^{N*}|}{N}\ge\alpha\bigg)\ge\mathbb{P}(Z=1).
\end{equation*}
When combined with (\ref{Z}) this completes the proof as we have exhausted
the possible values of $Z$.

We consider first the case $R\ge 1$. We can find
$t_0\in(1-\delta/2,1)$ such that $\beta(t_0)>\alpha +\delta/2$. Note
that $|X^N_{t_{0}}-x_{t_{0}}| \le \delta/2$ implies
\begin{align*}
|V^{N*}| /N&\ge X^{N,1}_{t_{0}}\ge t_0-\delta/2>1-\alpha,\\
|\Lambda^{N*}| /N&\ge X^{N,3}_{t_{0}}\ge\beta(t_0)-(1-t_0)\log(1-t_0)-\delta/2>\alpha.
\end{align*}
By Theorem \ref{t4}
\begin{equation*}
\liminf_{N\to\infty}\mathbb{P}\bigg(\sup_{t \le  t_0}
|X^N_t-x_t|\le\delta/2\bigg)\ge \mathbb{P}(Z>t_0)\ge\mathbb{P}(Z=1)
\end{equation*}
so we are done.

Consider next the case $R=0$. Fix $M\in\mathbb{N}$ and set
$\tilde\beta_j=\beta_j$ if $j \le  M$ and $\tilde\beta_j=0$
otherwise. Then, with obvious notation, we can choose $M$ so that
$\tilde{z}_0>1-\delta/2$, $\tilde\zeta=\emptyset$ and $\tilde{\beta}(z_0)>\alpha+\delta/2$.
Hence
\begin{equation*}
\mathbb{P}\bigg(\frac{|\tilde{V}^{N*}|}{N}\ge1-\delta\text{ and
}\frac{|\tilde{\Lambda}^{N*}|}{N}\ge\alpha\bigg)\to 1.
\end{equation*}
We can couple $\Lambda$ and $\tilde{\Lambda}$ so that
$\tilde{\Lambda}(A)=\Lambda(A)1_{|A|\le M}$. Then
$\tilde{V}^{N*}\subseteq V^{N*}$ and
$\tilde{\Lambda}^{N*}\le\Lambda^{N*}$, so this is enough.

There remains the case $R\in(0,1)$. In this case $\zeta$ is finite. We
have assumed that $R\notin\zeta$. So we can find $\rho\in(\sup\zeta,R)$
and $M\in{\mathbb N}$ such that, with obvious notation,
\begin{equation*}
\tilde{z}_0>1-\delta/2,\quad\tilde{\zeta}=\zeta,\quad\tilde{\beta}(\tilde{z}_0)>\alpha+\delta/2
\end{equation*} 
where $\tilde{\beta}(t), t\in[0,1)$, is defined by
\begin{equation*}
\tilde{\beta}(0)=\beta_0,\quad\tilde{\beta}'(0)= \beta_1,\quad
\tilde{\beta}''(t)=\begin{cases} \beta''(t), & t<\rho,\\
\sum^M_{j=2}j(j-1)\beta_j t^{j-2},& t\ge\rho.\end{cases}
\end{equation*}

Consider the collapsing hypergraph $(\tilde{\Lambda}^N_n)_{n\ge 0}$
which evolves as $(\Lambda^N_n)_{n\ge 0}$ up to $n=\nu(\rho)$, at
which time all hyperedges having at least two vertices and originally
having more than $M$ vertices are removed, so that $\tilde{\Lambda}^N_{\nu(\rho)}\le\Lambda^N_{\nu(\rho)}$.
After $\nu(\rho), (\tilde{\Lambda}^N_n)_{n\ge 0}$ evolves by selection
of patches as before. Denote by $\tilde{V}^{N*}$ the set of identifiable vertices 
 in $\tilde{\Lambda}^N_{\nu(\rho)}$ and by
$\tilde{\Lambda}^{N*}$ the corresponding identifiable hypergraph. Then 
\begin{equation*}
\tilde{V}^{N*}\subseteq V^{N*}\quad\hbox{and}\quad\tilde{\Lambda}^{N*}\le\Lambda^{N*}.
\end{equation*} 
\bigskip
A modification of Theorem \ref{t3} shows that
\begin{equation*}
{\mathbb P}\bigg(
|\tilde{X}^N_1-\tilde{x}_{\tilde{z}_0}|\le\delta/2\bigg)\quad\to\quad{\mathbb
P}(\tilde{Z}=\tilde{z}_0) ={\mathbb P}(Z=1)
\end{equation*}
with $\tilde{X}^N_1=(|\tilde{V}^{N*}| /N, 0, |\tilde{\Lambda}^{N*}|/N)$ 
and with
\begin{equation*}
\tilde{x}_t=(t,(1-t)(\tilde{\beta}(t)+\log(1-t)),\quad\tilde{\beta}(t)-(1-t)\log(1-t)).
\end{equation*}
All that changes in the proof is that, for $t\ge\rho$ the L\'evy kernel
is modified by replacing $\lambda_2$ by $\tilde{\lambda}_2$ given by
\begin{equation*}
\tilde{\lambda}_2(N,n) = N\sum^{n\wedge
(M-2)}_{i=0}\quad(i+1)(i+2)\beta_{i+2}\quad 
\frac{n(n-1)\dots(n-i+1)}{N(N-1)\dots(N-i+1)}.
\end{equation*}

The argument of Lemma 6.1 shows that for all $\rho'<1$ there is a
constant $C<\infty$ such that
\begin{equation*}
|N\tilde{\lambda}_2(N,n)-\tilde{\beta}''(n/N)|\quad \le \quad C/N
\end{equation*}
for all $N\in{\mathbb N}$ and $n=\lbrace 0,1,\dots,\lbrack N\rho'\rbrack\rbrace$.
Everything else is the same.
\bigskip
Now
$|\tilde{X}^N_1-\tilde{x}_{\tilde{z}_0}|\le\delta/2$ implies
\begin{equation*}
|V^{N*}|/N\ge |\tilde{V}^{N*}|/N
=\tilde{X}^{N,1}_1\ge\tilde{z}_0-\delta/2\ge 1-\delta,
\end{equation*}
\begin{equation*}
|\Lambda^{N*}|/N \ge |\tilde{\Lambda}^{N*}|/N = \tilde{X}^{N,3}_1
\ge\tilde{\beta}(\tilde{z}_0)-(1-\tilde{z}_0)\log (1-\tilde{z}_0)-\delta/2\ge\alpha
\end{equation*}
so
\begin{equation*}
\liminf_{N\to\infty}\quad{\mathbb P}\bigg(\frac{|V^{N*}|}{N} \ge
1-\delta\quad\hbox{and}\quad\frac{|\Lambda^{N*}|}{N}
\ge\alpha\bigg)\ge{\mathbb P}(Z=1)
\end{equation*}
as required.
\end{proof}
\begin{quote}
\emph{Acknowledgments}. The authors thank the 
Mathematisches Forschungsinstitut Oberwolfach for the invitation
to the meeting in 2000 at which this collaboration began, and
Brad Lackey for his interest in this project.
David Levin helped us to articulate some of the concepts presented herein.
\end{quote}
\end{document}